\documentclass[preprint]{amsart}
\usepackage{amsthm}
\RequirePackage{amssymb, amsfonts, amsmath, latexsym, verbatim, xspace, setspace, float}
\usepackage[toc,page]{appendix}

 \usepackage[all]{xy}
\usepackage{float}
\usepackage{tikz-cd}
\usepackage[outline]{contour}
\contourlength{4pt}
\usetikzlibrary{decorations.pathmorphing}
\tikzset{snake it/.style={decorate, decoration=snake}}
\usetikzlibrary{shapes, arrows}
\usepackage{mwe}

\usepackage{adjustbox}
\usetikzlibrary{arrows}
\usepackage{verbatim}
\usepackage{mathtools}
\usetikzlibrary{decorations.pathreplacing,decorations.markings}
\usepackage{ stmaryrd }
\usepackage{amsmath}
\usepackage{amssymb}
\usepackage{amsrefs}
\usepackage{hyperref}
\usetikzlibrary{intersections, matrix,decorations.pathreplacing}
\usetikzlibrary{calc}
\usepackage{longtable}
\usepackage{lipsum} 
\usepackage{mathtools}
\newtagform{black}[\color{black}]{(}{)}

\tikzset{
    ncbar angle/.initial=90,
    ncbar/.style={
        to path=(\tikztostart)
        -- ($(\tikztostart)!#1!\pgfkeysvalueof{/tikz/ncbar angle}:(\tikztotarget)$)
        -- ($(\tikztotarget)!($(\tikztostart)!#1!\pgfkeysvalueof{/tikz/ncbar angle}:(\tikztotarget)$)!\pgfkeysvalueof{/tikz/ncbar angle}:(\tikztostart)$)
        -- (\tikztotarget)
    },
    ncbar/.default=0.5cm,
}

\newcounter{markeq}
\newcommand{\pstrut}[1]{\vrule height0pt depth0pt width0pt #1 \fboxsep}
\newcommand*\bmarkeq{\stepcounter{markeq}%
  \tikz[remember picture]\node(startframe-\themarkeq){\pstrut{height}};%
  \kern\fboxsep}
\newcommand*\emarkeq{\kern\fboxsep
  \begin{tikzpicture}[remember picture,overlay]
    \node (endframe-\themarkeq){\pstrut{depth}};
    \draw[,red,opacity=0.8] (startframe-\themarkeq.north) 
      rectangle (endframe-\themarkeq.south);
  \end{tikzpicture}%
}

\usepackage{fp,amssymb}
\newsavebox\foobox

\usepackage{mwe}
\usepackage{tikz}

\tikzset{square left brace/.style={ncbar=0.5cm}}
\tikzset{square right brace/.style={ncbar=-0.5cm}}

\tikzset{round left paren/.style={ncbar=0.5cm,out=120,in=-120}}
\tikzset{round right paren/.style={ncbar=0.5cm,out=60,in=-60}}

\newtheorem{theorem}{Theorem}[section]
\newtheorem{lemma}[theorem]{Lemma}
\newtheorem{cor}[theorem]{Corollary}
\newtheorem{definition}[theorem]{Definition}
\newtheorem{prop}[theorem]{Proposition}

\theoremstyle{definition}

\newtheorem{example}[theorem]{Example}

\newtheorem{setting}[theorem]{Setting}

\newtheorem{convention}[theorem]{Convention}

\theoremstyle{remark}
\newtheorem{remark}[theorem]{Remark}

\numberwithin{equation}{section}

\newcommand{\partiald}[1] {\frac{\partial }{\partial #1}}

\DeclareFontFamily{U}{wncy}{}
    \DeclareFontShape{U}{wncy}{m}{n}{<->wncyr10}{}
    \DeclareSymbolFont{mcy}{U}{wncy}{m}{n}
    \DeclareMathSymbol{\Sh}{\mathord}{mcy}{"58} 

\newcommand\restr[2]{{
  \left.\kern-\nulldelimiterspace 
  #1 
  \vphantom{\big|} 
  \right|_{#2} 
  }}
\allowdisplaybreaks
\makeatletter
\newcommand{\leqnomode}{\tagsleft@true}
\newcommand{\reqnomode}{\tagsleft@false}
\makeatother

\begin{document}
 \pagestyle{headings}
  \title{The derivative of global surface-holonomy for a non-abelian gerbe}
  \author[C. Glass]{Cheyne J. Glass}
\address{St. Joseph's College New York, 155 W. Roe Blvd., Patchogue, NY 11772}
\email{cmiller5@sjcny.edu} 

 
 \begin{abstract}Starting with a non-abelian gerbe represented by a non-abelian differential cocycle, with values in a given crossed-module, this paper explicitly calculates a formula for the derivative of the associated surface holonomy of squares mapped into the base manifold; with spheres later considered as a special case. While the definitions in this paper used for gerbes, their connections, and the induced holonomy will initially be simplicial, translations into a cubical setting will be provided to aide in explicit coordinate-based calculations.  While there are many previously published results on the properties of these non-abelian gerbes, including some calculations of the derivative over a single open set, this paper endeavors to take these local calculations and glue them together across multiple open sets in order to obtain a single expression for the change in surface holonomy with respect to a one-parameter family of squares.  \end{abstract}
  \maketitle 

  \tableofcontents

\section{Introduction}
In \cite{TWZ}, an equivariantly closed differential form is associated to an abelian gerbe with connection by considering the derivative of the induced $2$-holonomy.  Originating as a first step in generalizing their work, this paper focuses on differentiating (Theorem \ref{dhol nice}) the global $2$-holonomy (Definition \ref{def of Hol}) for a non-abelian $\mathcal{G}$-gerbe (Definition \ref{def gerbe}).

Following Schreiber and Waldorf (\cite{SWI}, \cite{SWII}, \cite{SWIII}, and \cite{SWIV}), the local cocycle description for a non-abelian gerbe with connection on a smooth manifold is reviewed and adopted.  This paper uses their \emph{local transport data for bigons} implicitly in Section \ref{section square data}, but translates that data into \emph{local transport data for squares} as the process of differentiation became more manageable and organized in a cubical setting.  

The method for glueing this local data together to provide a global definition of $Hol$ (Definition \ref{def of Hol}) in a cubical setting was borrowed from Martins and Picken (specifically, Figure 3 in \cite{MP2}).  Their papers proved many properties for a global holonomy on the group-level which can be found in Section \ref{sec props of Hol}, where these properties are reviewed in addition to some other relevant observations.  This paper adds to those references by providing for a global formula, comprised of 
local data, for the derivative of global surface holonomy.

In Section \ref{sec dHol}, the main theorem of this paper, Theorem \ref{dhol nice}, is stated and proven, which essentially reads 
$$d(Hol) = Hol \cdot \int_{Sq} H \quad \text{(modulo terms on the boundary of $Sq$)},$$
where $\int_{Sq} H$ represents fiber-integration of the $3$-curvature terms for the given non-abelian gerbe through the interior of the square.  The proof of this theorem, given in Section \ref{dHol nice proof}, amounts to considering a $1$-parameter family of squares, considering the associated arrangement of cubes from the local data, and organizing the terms in the derivative accordingly.  

Finally, in Sections \ref{sec spheres} and \ref{sec abelian}, some examples are offered where Theorem \ref{dhol nice} has an even cleaner representation.  In the case where the surfaces which are integrated over are spheres, $d(Hol)$ has only a boundary term at the base point.  In the case where the gerbe is abelian, the well known situation where there are no terms for $d(Hol)$ on the boundary is reprodouced.  

The hope after this paper is to continue the work of finding a non-abelian analogue of the work done in \cite{TWZ}, via a subsequent paper which will use the result of Theorem \ref{dhol nice}, and the appropriate cohomology theory, to find some equivariantly closed element representing a non-abelian gerbe with connection.  In a separate project, the goal is to extend this paper's derivative for $2$-dimensional holonomy (landing in crossed-modules) to the derivative of $3$-dimensional holonomy (landing in the appropriate version of a $2$-crossed module).  Furthermore, some of the geometric comments regarding the ``globalness'' of this 2-holonomy are planned to be formalized in current joint work with Micah Miller, Thomas Tradler, and Mahmoud Zeinalian.
\subsection*{Acknowledgement}
The author would like to thank Thomas Tradler for many helpful conversations about this paper.  

\section{Conventions, Notation, and Setup}\label{section setup}
In order to arrive at a definition for global 2-holonomy, Definition \ref{def of Hol}, it is necessary to introduce some preliminary definitions and conventions.
\subsection{Diffeological Spaces}\label{sec diffeo}
Since we will be working on the space of smooth maps, $M^{Sq}:= \{ \Sigma: Sq \to M\}$, where $Sq:= [0,1] \times [0,1]$ is the \emph{standard square} and $M$ is a smooth manifold, it is convenient to use the language of \emph{diffeological spaces} and \emph{plots} as described in \cite{BaHo} and \cite{IZ} (see \cite{C2} for an early reference on ``plots'').  For our purposes, $M^{Sq}$ is a \emph{diffeology} whose \emph{plots} are given by  maps $P : U \to C^{\infty}(M, Sq)$ which are smooth in both variables, i.e. the map $(r,x) \mapsto P(r)(x)$ is a smooth map from the subset $U \subset \mathbb{R}^n$ to $M$.

\subsubsection{Covering $M^{Sq}$ With Open Sets $\mathcal{N}$}

It is possible that any square mapped into $M$ will not be entirely contained within one open set, $U_{\alpha}$, of a chosen cover, $\mathcal{U}$.  Instead, the square, $Sq$, is subdivided into a grid whereby each sub-square, $Sq_{(p,q)}$, is mapped entirely inside some open set $U_{i_{(p,q)}} \in \mathcal{U}$, as was the key idea borrowed from \cite{MP2} and is now explained below.  Note that while a finite collection of open sets is usually indexed, $U_{i_1}, U_{i_2},\ldots,  U_{i_n}$, it is helpful in this context to use the indexing convention, $i_{(p,q)}$, which will range over the grid: $U_{i_{(1,1)}}$, $U_{i_{(1,2)}}$, $\ldots, U_{i_{(2,1)}}$, $U_{i_{(2,2)}}$, $\ldots, U_{i_{(n,m)}}$.
\begin{definition}\label{def of N}
For a fixed open cover, $\mathcal{U}$, of $M$, a choice of grid $I = \{1, \ldots , n\} \times \{1, \ldots, m\}$, and a choice of open sets $\{U_{i_{(p,q)}}\}_{(p,q) \in I}$, define
$$\mathcal{N}_I := \{\Sigma \in M^{Sq} \mid \ \text{ for each } (p,q) \in I, \Sigma(Sq_{(p,q)}) \subset U_{i_{(p,q)}}\}.$$ 
where $Sq_{(p,q)} := \left[\frac{p-1}{n}, \frac{p}{n}\right] \times \left[ \frac{q-1}{m}, \frac{q}{m} \right] \subset [0,1] \times [0,1]$.  When the indexing set, $I$, is understood, we will simply refer to the open set as $\mathcal{N}$.
\end{definition}
It will be useful later in Section \ref{section glue Hols} to have the notion of a \emph{grid} on $\Sigma \in \mathcal{N}_I$.
\begin{definition}\label{def of grid}
Given a square, $\Sigma \in \mathcal{N}_I$, define the \emph{grid on} $\Sigma$ by the following set of data: For each $i= (p, q) \in I$ define the $i$-\underline{face} by $\Sigma_i:= \restr{\Sigma}{Sq_i}$; For each $i = (p, q), j=(p +1, q)$ define the $ij$-\underline{vertical edge}, $\gamma^v_{ij}$, by 
$$\gamma^v_{ij}:= \restr{\Sigma_i}{\{\frac{p}{n} \} \times \left[ \frac{q -1}{m}, \frac{q}{m}\right] }=  \restr{\Sigma_j}{\{\frac{p}{n}\} \times \left[ \frac{q -1}{m}, \frac{q}{m}\right]} =  \restr{\Sigma}{\{ \frac{p}{n}\} \times \left[ \frac{q-1}{m}, \frac{q}{m}\right]};$$
For each $i = (p, q), j=(p, q+1)$ define the $ij$-\underline{horizontal edge}, $\gamma^h_{ij}$, by 
$$\gamma^h_{ij}:=\restr{\Sigma_i}{\left[ \frac{p -1}{n}, \frac{p}{n}\right]  \times \{\frac{q}{m}\}}  =\restr{\Sigma_j}{\left[ \frac{p -1}{n}, \frac{p}{n}\right]  \times \{\frac{q}{m}\}}  =  \restr{\Sigma}{\left[ \frac{p -1}{n}, \frac{p}{n}\right]  \times \{\frac{q}{m}\}};$$
Any such $ij$-edge has a source vertex and target vertex labeled by $x_{ij}^0$ and $x_{ij}^1$, respectively.  For each $ijkl$ where $i=(p, q), j=(p, q +1), k=(p +1, q), l=(p +1, q + 1)$ define the $ijkl$- \underline{vertex}, $x_{ijkl}$, by 
$$x_{ijkl}:= \restr{\Sigma}{(\frac{p}{n}, \frac{q}{m})}.$$
For each $i_N = (p,1)$, $i_S=(p,m)$,  $i_W=(1,q)$, and $i_E=(n,q)$, define the \underline{northern boundary edge}, $\gamma_i^N$, the  \underline{southern boundary edge}, $\gamma_i^S$, the \underline{western boundary edge}, $\gamma_i^W$, and the \underline{eastern boundary edge}, $\gamma_i^E$, respectively by:
$$\gamma^N_{i}:=\restr{\Sigma_i}{\left[ \frac{p -1}{n}, \frac{p}{n}\right]  \times \{ 0 \}} \quad \gamma^S_{i}:=\restr{\Sigma_i}{\left[ \frac{p -1}{n}, \frac{p}{n}\right]  \times \{1 \}} \quad \gamma^W_{i}:=\restr{\Sigma_i}{  \{0 \} \times \left[ \frac{q -1}{m}, \frac{q}{m}\right]} \quad \gamma^E_{i}:=\restr{\Sigma_i}{  \{1 \} \times \left[ \frac{q -1}{m}, \frac{q}{m}\right]}.  $$
\end{definition}

Every diffeology on $X$ induces a topology on $X$, called the $\mathcal{D}$-topology.  Page 54 of \cite{IZ} states the following characterization of the $\mathcal{D}$-topology:
\begin{prop} A subset $A$ of a diffeological space, $X$, is open for the $\mathcal{D}$-topology if and only if for every plot, $\rho: U \to X$, $\rho^{-1}(A)$ is open in $U$.\end{prop}
It is straightforward to check that the sets $\mathcal{N}$ are open in the diffeology $M^{Sq}$ under its $\mathcal{D}$-topology:
\begin{prop}
Consider the diffeological space $M^Y: = \{ f : Y \to M \mid f \text{ is smooth }\}$, with plots $\rho: U \to M^Y$ given by smooth maps $\tilde{\rho}: U \times Y \to M$, where $(r,y) \mapsto \rho(r)(y).$  If $K \subset Y$ is compact and $U \subset M$ is open, then $$\mathcal{N} (K,U):= \{f \in M^Y \mid \restr{f}{K} \subset U \}$$ is an open set in $M^Y$.
\end{prop}
\begin{cor}
Each $\mathcal{N}_I$ is an open subset of $M^{Sq}$.
\end{cor}

\begin{prop}\label{prop Ns cover MSq}
For a fixed open cover $\mathcal{U}$, the open sets $\mathcal{N}_I$, ranging over all choices from Definition \ref{def of N}, cover $M^{Sq}$.
\end{prop}

\subsection{Crossed Module Conventions and Relations}
An early reference for the understanding that crossed modules were helpful structures for dealing with 2-dimensional algebra is \cite{BrSp}.  In this section, following \cite{GiPf} as a reference for notation and convention, the definition of crossed modules, and some related properties that will be essential later on, are now reviewed.

\begin{definition}\label{CM LG}
A crossed module of Lie Groups is a pair of Lie groups, $(H,G)$, with a smooth group homomorphism, $(H \xrightarrow{t} G)$, called the target, and an action $\alpha: G \to Aut(H)$, written $\alpha_g(h)$, so that $t$ and $\alpha$ are required to satisfy the following compatibility coniditions
\begin{align}
t(\alpha_g(h)) &= g t(h) g^{-1}\\
\alpha_{t(h)}(h') &= h h' h^{-1} \label{eq alpha t h commutes}
\end{align}
\end{definition}
\begin{example}\label{ex AutH}
If $H$ is a Lie group, then $G:= Aut(H)$ induces a crossed module $(H \xrightarrow{t} Aut(H))$ via the target $t: H \to Aut(H)$ given by $$t(h)(h'):= hh'h^{-1}$$ and the action $\alpha: G \to Aut(H)$ given by the identity automorphism.
\end{example}

\begin{example}\label{ex BS1}
Define the crossed module $\mathcal{B}S^{1}:= (S^1 \xrightarrow{t} \{ *\})$ given by the trivial target and identity action $\alpha$.  This is the crossed module often used in the study of abelian gerbes.
\end{example}

Associated to such a crossed module is a crossed module of Lie algebras:
\begin{definition}\label{CM LA}
A crossed module of Lie algebras is a pair of Lie algebras, $(\mathfrak{h},\mathfrak{g})$, with a Lie algebra map, $(\mathfrak{h} \xrightarrow{t} \mathfrak{g})$, called the target, and a map $\alpha: \mathfrak{g} \to \partial er(\mathfrak{h})$, written $\alpha_A(B)$, so that the two maps $t$ and $\alpha$ must satisfy the following compatibility conditions:
\begin{align}
t(\alpha_XY) &= [X, t(Y)]\\
\alpha_{t(Y_1)}(Y_2) &= [Y_1, Y_2]
\end{align}
\end{definition}
A useful proposition regarding the center of $H$ is as follows:
\begin{prop}\label{trivial t center}
The kernel of the target map, $t: H \to G$, is in the center, $Z(H)$, of $H$.  Similarly, the kernel of the target map, $t: \mathfrak{h} \to \mathfrak{g}$ is in the cetner, $Z(\mathfrak{h})$, of $\mathfrak{h}$.
\begin{proof}
If $t(h) = 1$, then for any $h' \in H$, $h' = \alpha_1(h') = \alpha_{t(h)}(h')= h h' h^{-1}$, thus $h \in Z(H)$.  A similar proof can be applied to the statement for $X \in \mathfrak{h}$ with $t(X) = 0$:
$$0 = \alpha_0(Y) = \alpha_{t(X)}(Y) = [X,Y].$$
\end{proof}
\end{prop}

An important and well-known lemma, \cite[Proposition~1.4]{KN1} in calculating the derivative of the $\alpha$ map is provided now for later reference
\begin{lemma}\label{lemma d alpha}
For the function $\alpha: G \times H \to H$ defined above, given functions $g(t): \mathbb{R} \to G$ and $h(t): \mathbb{R} \to H$, we have 
\begin{equation}
\restr{\partiald{t}}{t=t_0}\left(\alpha_{g(t)}(h(t))\right) = (\alpha_{g(t_0)})_*\left( \restr{\partiald{t}}{t=t_0}h(t)\right) + (\alpha_{h(t_0)})_*\left( \restr{\partiald{t}}{t=t_0}g(t)\right)
\end{equation} as an equality in $T_{\alpha_{g(t_0)}(h(t_0))}H$, the tangent space of $H$ at $\alpha_{g(t_0)}(h(t_0))$.  We will sometimes write 
$$\alpha_{\left(\restr{\partiald{t}}{t=t_0}g(t)\right) }h(t_0) := (\alpha_{h(t_0)})_*\left( \restr{\partiald{t}}{t=t_0}g(t)\right)$$
and refer to these as ``\emph{path terms}''.
\end{lemma}

\subsection{Local Differential Data for a Gerbe}\label{section local cocycle data}
In \cite{NW} it is shown that the cocycle description of a non-abelian gerbe is equivalent to the other three common formulations: classifying maps, groupoid bundle gerbes, and principal 2-bundles.  In \cite{SWIII} they go on to provide their simplicial formulation for a connection on, and associated parallel transport for, a given non-abelian gerbe.  Following \cite{SWIII}, and since it is assumed we are working on a \emph{good open cover}, i.e. one where all open sets and their $n$-fold intersections are contractible, the following (normalized) local data for a gerbe are used:
\begin{definition}\label{def gerbe}
Given a smooth manifold, $M$, an open cover $\mathcal{U} = \{ U_i\}$ of $M$, and a crossed module of Lie groups $\mathcal{G} = (H \xrightarrow{t} G)$, a $\mathcal{G}$-gerbe with a connection on $M$, subordinate to the cover $\mathcal{U}$, is defined by the following (normalized) local cocycle data:
\begin{itemize}
\item {\bf On each open set}, $U_i$ a pair $$(A_i \in \Omega^1(U_i, \mathfrak{g}), B_i \in \Omega^2(U_i, \mathfrak{h})).$$
\item {\bf On each intersection}, $U_{ij}:= U_i \cap U_j$ a pair $$(g_{ij} \in \Omega^0(U_{ij}, G), a_{ij}\in \Omega^1(U_{ij}, \mathfrak{h})).$$
\item {\bf On each triple intersection}, $U_{ijk}:= U_i \cap U_j \cap U_k$ a function $$f_{ijk} \in \Omega^0(U_{ijk}, H).$$
\end{itemize}
satisfying the relations,
\begin{enumerate}
\item {\bf On each open set}, \begin{align*}
R_i:= dA_i + \frac{1}{2}[A_i \wedge A_i] &= t(B_i)\\
g_{ii} &=1\\
a_{ii} &=0.
\end{align*}
\item {\bf On each intersection}, \begin{align*}
A_j &=g_{ij} A_i g_{ij}^{-1} - dg_{ij} g_{ij}^{-1} - t(a_{ij})\\
B_j &= \alpha_{g_{ij}}(B_i) - \nabla_j(a_{ij}) - \frac{1}{2}[a_{ij} \wedge a_{ij}]\\
f_{iij} &= f_{ijj} = 1.
\end{align*}
\item {\bf On each triple intersection}, \begin{align*}
g_{ik} &= t(f_{ijk}) g_{jk}g_{ij}\\
f_{ijk}\cdot a_{ik}\cdot f_{ijk}^{-1} &= (\alpha_{g_{jk}})_*(a_{ij}) + a_{jk} + ((\alpha_{f_{ijk}})_*(A_k))\cdot f_{ijk}^{-1} + df_{ijk}\cdot f_{ijk}^{-1}.
\end{align*}
\item {\bf On each quadruple intersection}, \begin{align*}
f_{ikl}\alpha_{g_{kl}}(f_{ijk})&= f_{ijl}f_{jkl}.\\
\end{align*}
where $\nabla_j(\omega):= d\omega + \alpha_{A_j}( \omega)$,  $g\cdot A:= (L_g)_*(A)$, and $A \cdot g:= (R_g)_*(A)$.
\end{enumerate}
\end{definition}
Girelli and Pfeiffer \cite{GiPf}, as well as  Baez and Schreiber \cite{BaSc}, define the curvature 3-form of a gerbe $H_i \in \Omega^3(U_i, \mathfrak{h})$ by 
$$H_i:= \nabla_i(B_i):= dB_i + \alpha_{A_i}(B_i)$$ where $B_i \in \Omega^2(U_i, \mathfrak{h})$ and $(A_i,B_i)$ define the 2-connection of our gerbe on $U_i$.  It is straight-forward to check the following  properties of the curvature $3$-form.
\begin{prop}\label{props of 3 curvature}
The curvature 3-form, $H_i$, of a gerbe, satisfies
\begin{enumerate}
\item $t(H_i) = 0$
\item $\nabla_i(H_i) = 0$.
\item $H_j = \alpha_{g_{ij}}(H_i)$, on each intersection of open sets, $U_{ij}$.
\end{enumerate}
\end{prop}

\subsection{Local 2-Holonomy}
In this section, the simplicial definition of gerbes will give rise to a cubical notion of 2-holonomy.  

\subsubsection{Comments on the local geometry of 2-holonomy}
In the works of Schreiber and Waldorf \cite{SWI, SWII, SWIII} as well as Martins and Picken \cite{MP1, MP2}, the higher-holonomy is constructed as a map of 2-groupoids from a certain path 2-groupoid (see \cite[section 2.1]{SWII} or \cite[section 2.3]{MP2}) to the 2-groupoid described by a crossed module of Lie groups with on object (see \cite[Theorem A]{BrSp} for this second 2-groupoid).  Having such a map requires certain coherence conditions: for example that the algebraic target for the 2-holonomy of a surface be related to the 1-holonomy of the geometric boundary for that surface; and that the holonomy of a composition of surfaces/paths corresponds to the algebraic product of the holonomies of these surfaces/paths.  

In both of these sets of authors' respective path 2-groupoids (Martins and Picken using a cubical special double groupoid a la Brown and Spencer versus Schreiber and Waldorf using a simplicial double groupoid), the objects are points in the base manifold, and the 1-morphisms are \emph{thin homotopy classes} of paths.  That is, the morphisms are equivalence classes of paths where the relation is given by certain homotopies between two paths having at most rank 1 \cite[Definition 1.1]{SWII}.  Since this paper formally uses the data from \cite{SWII}, their simplicial path 2-groupoid is technically being used.  However, by considering a square mapped into a manifold as a bigon whose source is the $(t=0, s=0)$-vertex of the square and whose target is the boundary of the square, we can apply the higher holonomy maps of \cite{SWII} to our squares mapped into the manifold and proceed from there.

\subsubsection{Local Transport Data for Squares}\label{section square data}
Fix a gerbe on $M$ subordinate to the cover $\mathcal{U}$, as defined in Definition \ref{def gerbe}.  Then, following a modified version of the transport data for bigons given in \cite[Section 1.3]{SWII}, define smooth functions $$hol_i: U_i^{[0,1]} \to G, \quad Hol_i: U_i^{Sq} \to H, \quad Hol_{ij}: U_{ij}^{[0,1]} \to H, \text{ and } \quad Hol_{ijkl}: U_{ijkl} \to H.$$   These local functions satisfy the differential equations below.  Essential to the construction of these functions in \cite{SWII} are the concepts that (a) each path, $\gamma: [0,1] \to U_i$, can be extended trivially to a path with the same image in $M$, $\tilde{\gamma}$, having domain $\mathbb{R}$, and (b) from such a path, $\tilde{\gamma}: \mathbb{R} \to U_i$, we can consider the one-parameter family of paths, $\tilde{\gamma}_t$, where for each $t \in \mathbb{R}$, the path's image is that of $\restr{\gamma}{[0,t]}$ while being re-parametrized to still satisfy the definition of being a morphism in the \emph{path groupoid} of $M$ \cite[Section 1.1]{SWII}.  Similarly, squares $\Sigma \in U_i^{Sq}$ produce one-parameter families of squares $\Sigma_s$ whose image is that of $\Sigma_{[0,1] \times [0,s]}$ while being 2-morphisms in the \emph{path 2-groupoid} of $M$ \cite[Section 2.1]{SWII}.
\begin{itemize}
\item[(D1)] For a one-parameter family of paths, $\gamma_t$, induced by the path, $\gamma(t)$:
$$\partiald{t}hol_i(\gamma_t) ^{-1}= hol_i \cdot A_i \left(\partiald{t}\right)$$
\item[(D2)]\label{Holi diff eq} For a one-parameter family of squares, $\Sigma_{s}$, induced by the square, $\Sigma(t,s)$:
$$\partiald{s}Hol_i(\Sigma_{s}) = Hol_i(\Sigma_s) \cdot \int_0^1( \alpha_{hol_i(-,s)})_*(B_i) dt \left(\partiald{s}\right)$$
\item[(D3)] \label{derivative of Holij at 0} For a one-parameter family of paths, $\gamma_t$, induced by the path, $\gamma(t)$:
$$\restr{\partiald{t}}{t=0} (Hol_{ij}(\gamma_t)) = \alpha_{g_{ij}^{-1}(\gamma(0))}\left(\restr{a_{ij}}{\gamma(0)} \left(\restr{\partiald{t}}{t=0} \right) \right)$$
\end{itemize}
The collection of local transport data for squares $Hol_i$ and $Hol_{ij}$ have the following targets:
\begin{itemize}
\item For $Hol_i$:
$$t(Hol_i) = hol_i(-,0)^{-1}hol_i(t,-)^{-1}hol_i(-,s)hol_i(0,-)$$
 which will be pictured diagrammatically as:
\begin{equation}\label{target Holi}
\begin{tikzpicture}[baseline={([yshift=-.5ex]current bounding box.center)},vertex/.style={anchor=base,
     circle,fill=black!25,minimum size=18pt,inner sep=2pt}]   
\fill (0,3) circle (2pt);
\path[thick] (0,3)       edge coordinate[near end] (h)node[above]  {\contour{white}{$hol_i(\gamma(-,0))^{-1}$}}(3,3);
\path[thick]  (3,3)       edge coordinate[near end] (h)node[right]  {\contour{white}{$hol_i(\gamma(t,-))^{-1}$}}(3,0);
\path[thick] (0,0)       edge coordinate[near end] (h)node[below] {\contour{white}{$hol_i(\gamma(-,s))$}}(3,0);
\path[thick, shorten <=0.5cm,<-]    (0,3)       edge coordinate[near end] (h)node[left]  {\contour{white}{$hol_i(\gamma(0,-))$}}(0,0);
\node at (1.5, 1.5) {$Hol_i$};
\end{tikzpicture}
\end{equation}
\item For $Hol_{ij}$:
$$t(Hol_{ij})  = hol_i^{-1} g_{ij}^{-1}(y) hol_j g_{ij}(x)$$
which will be pictured diagrammatically as:
\begin{equation}\label{target Holij}
\begin{tikzpicture}[baseline={([yshift=-.5ex]current bounding box.center)},vertex/.style={anchor=base,
     circle,fill=black!25,minimum size=18pt,inner sep=2pt}]   
\fill (0,3) circle (2pt);
\path[thick] (0,3)       edge coordinate[near end] (h)node[above]  {\contour{white}{$hol_i(\gamma)^{-1}$}}(3,3);
\path[thick]  (3,3)       edge coordinate[near end] (h)node[right]  {\contour{white}{$g_{ij}(\gamma(1))^{-1}$}}(3,0);
\path[thick] (0,0)       edge coordinate[near end] (h)node[below] {\contour{white}{$hol_j(\gamma)$}}(3,0);
\path[thick, shorten <=0.5cm,<-]    (0,3)       edge coordinate[near end] (h)node[left]  {\contour{white}{$g_{ij}(\gamma(0))$}}(0,0);
\node at (1.5, 1.5) {$Hol_{ij}$};
\end{tikzpicture}
\end{equation}
\item For $Hol_{ijkl}:= \alpha_{g_{ik}^{-1}}\left(\alpha_{g_{kl}^{-1}}(f_{jkl}) f_{ijk}^{-1} \right)$:
$$t(Hol_{ijkl}) = g^{-1}_{ik}g^{-1}_{kl}g_{jl}g_{ij}$$
 which will be pictured diagrammatically as:
\begin{equation}\label{target Holijkl}
\begin{tikzpicture}[baseline={([yshift=-.5ex]current bounding box.center)},vertex/.style={anchor=base,
     circle,fill=black!25,minimum size=18pt,inner sep=2pt}]   
\fill (0,3) circle (2pt);
\path[thick] (0,3)       edge coordinate[near end] (h)node[above]  {\contour{white}{$g^{-1}_{ik}$}}(3,3);
\path[thick]  (3,3)       edge coordinate[near end] (h)node[right]  {\contour{white}{$g^{-1}_{kl}$}}(3,0);
\path[thick] (0,0)       edge coordinate[near end] (h)node[below] {\contour{white}{$g_{jl}$}}(3,0);
\path[thick, shorten <=0.5cm,<-]    (0,3)       edge coordinate[near end] (h)node[left]  {\contour{white}{$g_{il}$}}(0,0);
\node at (1.5, 1.5) {$Hol_{ijkl}$};
\end{tikzpicture}
\end{equation}
\end{itemize}

\subsubsection{Glueing-Paths}
Some of the conventions and computations which are necessary to glue the local transport data for squares defined above are laid in this section.  Much of the details in this section can be described via the following proposition which can be seen as a result of the works of Brown, Spencer, and Higgins in \cite{BrSp, BrHi, Br}.
\begin{prop}[{\cite[Theorem A]{BrSp}}]
Let $\mathcal{D}\mathcal{G}$ be the category whose objects are special double groupoids with special connection and whose arrows are the morphisms of double groupoids preserving the connection, and $\mathcal{C}$ be the category of crossed modules.  If we consider the full sub-category of $\mathcal{D}\mathcal{G}$ whose double groupoids have exactly one object labeled, $\mathcal{D}\mathcal{G}^!$, then there is an equivalence of categories $\gamma:  \mathcal{D}\mathcal{G}^! \to \mathcal{C}$.
\end{prop}
As a medium of discussion, the interchange law for our ``glueing paths'' is proved.  At the end of this section, the established conventions are used to prove Propositions \ref{Vertex Cube} and \ref{Edge Cube}.

If two (horizontally) adjacent squares are to glued, define horizontal multiplication by using the following procedure\footnote{This procedure is helpful to be conceptualized as a ``zip/unzip'' procedure, in order to easily follow some calculations.}:  
\begin{align}\label{horizontal glue}
\begin{tikzpicture}[baseline={([yshift=-.5ex]current bounding box.center)},vertex/.style={anchor=base,
    circle,fill=black!25,minimum size=18pt,inner sep=2pt}]
\fill (0,2.2) circle (2pt);
\draw[dotted] (0,2) to coordinate[midway] (a) (2,2) to coordinate[midway] (b) (2,0) to coordinate[midway] (c) (0,0) to coordinate[midway] (d) (0,2);
\draw[dotted] (2,2) to  coordinate[midway] (e) (4,2) to  coordinate[midway] (f) (4,0) to  coordinate[midway] (g) (2,0) to (2,2);
\node at (1,1) {$Hol_1$};
\node at (3,1) {$Hol_2$};
\draw[shorten >=0.2cm, ->] (0,2.2) to (2,2.2) to (2,0.2) to (0, 0.2)  to (0,2) to (-0.2,2) to (-0.2,0) to  (2.2,0) to (2.2,2) to (4,2) to (4,-0.2) to (-0.4,-0.2) to (-0.4,2);
\node [fill=white] at (a) {$a$};
\node [fill=white] at (b) {$b$};
\node [fill=white] at (c) {$c$};
\node [fill=white] at (d) {$d$};
\node [fill=white] at (e) {$e$};
\node [fill=white] at (f) {$f$};
\node [fill=white] at (g) {$g$};
\end{tikzpicture}
=   \begin{tikzpicture}[baseline={([yshift=-.5ex]current bounding box.center)},vertex/.style={anchor=base,
    circle,fill=black!25,minimum size=18pt,inner sep=2pt}]
\draw[shorten >=0.2cm, ->] (0,2) to coordinate[midway] (a) (3,2) to coordinate[midway] (e) (6,2) to coordinate[midway] (f) (6, 0)  to coordinate[midway] (g) (3,0) to coordinate[midway]  (c) (0,0) to coordinate[midway]  (d) (0,2);
\draw[-implies,double equal sign distance, shorten >= 0.2cm, shorten <= 0.2cm] (3,2) -- node[fill=white] {$Hol_1\cdot \alpha(d^{-1}c^{-1}b)(Hol_2)$} (3,0) ;
\node [above] at (a) {$a$};
\node [above] at (e) {$e$};
\node [below] at (c) {$c$};
\node [left] at (d) {$d$};
\node [below] at (g) {$g$};
\node [right] at (f) {$f$};
  \end{tikzpicture}
\end{align}
Similarly, for vertically adjacent squares, vertical multiplication is defined by using the following procedure:
\begin{align}\label{vertical glue}
\begin{tikzpicture}[baseline={([yshift=-.5ex]current bounding box.center)},vertex/.style={anchor=base,
    circle,fill=black!25,minimum size=18pt,inner sep=2pt}]
\fill (0,4) circle (2pt);
\draw[white] (0,4) to coordinate[midway] (d) (0,2) to coordinate[midway] (c) (2,2) to coordinate[midway] (b) (2,4) to coordinate[midway] (a) (0,4);
\draw[white] (0,2) to  coordinate[midway] (h) (0,0) to  coordinate[midway] (g) (2,0) to  coordinate[midway] (f) (2,2) to (0,2);
\node at (1,1) {$Hol_2$};
\node at (1,3) {$Hol_1$};
\draw[shorten >=0.2cm, ->] (0.2,4) to (2,4) to (2,2.2) to (0.2, 2.2)  to (0.2,3.8) to (0,3.8) to (0,1.8) to  (2,1.8) to (2,0) to (-0.2,0) to (-0.2,4);
\node [fill=white] at (a) {$a$};
\node [fill=white] at (b) {$b$};
\node [fill=white] at (c) {$c$};
\node [fill=white] at (d) {$d$};
\node [fill=white] at (f) {$f$};
\node [fill=white] at (g) {$g$};
\node [fill=white] at (h) {$h$};
\end{tikzpicture}
=   \begin{tikzpicture}[baseline={([yshift=-.5ex]current bounding box.center)},vertex/.style={anchor=base,
    circle,fill=black!25,minimum size=18pt,inner sep=2pt}]
\fill (0,4) circle (2pt);
\draw[white] (0,4) to coordinate[midway] (d) (0,2) to coordinate[midway] (c) (4,2) to coordinate[midway] (b) (4,4) to coordinate[midway] (a) (0,4);
\draw[white] (0,2) to  coordinate[midway] (h) (0,0) to  coordinate[midway] (g) (4,0) to  coordinate[midway] (f) (4,2) to (0,2);
\node at (2,2) {$Hol_1 \cdot \alpha(a^{-1})(Hol_2)$};
\draw[shorten >=0.2cm, ->] (0,4) to (4,4) to (4,0) to (0, 0)  to (0,4);
\node [fill=white] at (a) {$a$};
\node [fill=white] at (b) {$b$};
\node [fill=white] at (d) {$d$};
\node [fill=white] at (h) {$h$};
\node [fill=white] at (f) {$f$};
\node [fill=white] at (g) {$g$};
\end{tikzpicture}
\end{align}

\begin{prop}\label{prop interchange}
The procedures above for composing 2-squares satisfy the ``\emph{interchange law}''.  Moreover, any grid of squares glues together to a single square providing a unique element in $H$ associated to the grid.
\begin{proof}
The proof is exhibited in diagram-form.  When first composing horizontally and then composing vertically:
\begin{align*}
&\begin{tikzpicture}[ scale=1.2,baseline={([yshift=-.5ex]current bounding box.center)},vertex/.style={anchor=base,
    circle,fill=black!25,minimum size=18pt,inner sep=2pt}]
\fill (0,4) circle (2pt);
\draw[shorten <= 0.2cm, shorten >= 0.2cm, ->] (0,4) to coordinate[midway]  (a) (2,4);
\draw[shorten <= 0.2cm, shorten >= 0.2cm, ->] (2,4) to coordinate[midway]  (b) (4,4);
\draw[shorten <= 0.2cm, shorten >= 0.2cm, ->] (0,4) to coordinate[midway]  (c) (0,2);
\draw[shorten <= 0.2cm, shorten >= 0.2cm, ->] (2,4) to coordinate[midway]  (d) (2,2);
\draw[shorten <= 0.2cm, shorten >= 0.2cm, ->] (4,4) to coordinate[midway]  (e) (4,2);
\draw[shorten <= 0.2cm, shorten >= 0.2cm, ->] (0,2) to coordinate[midway]  (f) (2,2);
\draw[shorten <= 0.2cm, shorten >= 0.2cm, ->] (2,2) to coordinate[midway]  (g) (4,2);
\draw[shorten <= 0.2cm, shorten >= 0.2cm, ->] (0,2) to coordinate[midway]  (h) (0,0);
\draw[shorten <= 0.2cm, shorten >= 0.2cm, ->] (2,2) to coordinate[midway]  (i) (2,0);
\draw[shorten <= 0.2cm, shorten >= 0.2cm, ->] (4,2) to coordinate[midway]  (j) (4,0);
\draw[shorten <= 0.2cm, shorten >= 0.2cm, ->] (0,0) to coordinate[midway]  (k) (2,0);
\draw[shorten <= 0.2cm, shorten >= 0.2cm, ->] (2,0) to coordinate[midway]  (l) (4,0);
 \foreach \x in {a,b,c,...,l} {\node [fill=white] at (\x) {$\x$};}
 \node at (1,3) {$Hol_1$};
  \node at (3,3) {$Hol_2$};
   \node at (1,1) {$Hol_3$};
    \node at (3,1) {$Hol_4$};
\end{tikzpicture}
= \begin{tikzpicture}[ scale=1.2, baseline={([yshift=-.5ex]current bounding box.center)},vertex/.style={anchor=base,
    circle,fill=black!25,minimum size=18pt,inner sep=2pt}]
\fill (0,4) circle (2pt);
\draw[shorten <= 0.2cm, shorten >= 0.2cm, ->] (0,4) to coordinate[midway]  (ba) (6,4);
\draw[shorten <= 0.2cm, shorten >= 0.2cm, ->] (0,4) to coordinate[midway]  (c) (0,2);
\draw[shorten <= 0.2cm, shorten >= 0.2cm, ->] (6,4) to coordinate[midway]  (e) (6,2);
\draw[shorten <= 0.2cm, shorten >= 0.2cm, ->] (0,2) to coordinate[midway]  (gf) (6,2);
\draw[shorten <= 0.2cm, shorten >= 0.2cm, ->] (0,2) to coordinate[midway]  (h) (0,0);
\draw[shorten <= 0.2cm, shorten >= 0.2cm, ->] (6,2) to coordinate[midway]  (j) (6,0);
\draw[shorten <= 0.2cm, shorten >= 0.2cm, ->] (0,0) to coordinate[midway]  (lk) (6,0);
 \foreach \x in {ba, c, e, gf, h, j, lk} {\node [fill=white] at (\x) {$\x$};}
 \node at (3,3) {$Hol_1 \cdot \alpha(c^{-1}f^{-1}d)(Hol_2)$};
  \node at (3,1) {$Hol_3 \cdot \alpha(h^{-1}k^{-1}i)(Hol_4)$};
\end{tikzpicture}\\
= & \begin{tikzpicture}[ scale=1.2, baseline={([yshift=-.5ex]current bounding box.center)},vertex/.style={anchor=base,
    circle,fill=black!25,minimum size=18pt,inner sep=2pt}]
\fill (0,2) circle (2pt);
\draw[shorten >= 0.2cm, ->] (0,2) to coordinate[midway]  (ba) (11,2) to coordinate[midway] (je) (11,0) to coordinate[midway] (lk) (0,0) to coordinate[midway] (hc) (0,2);
 \foreach \x in {ba, je, lk, hc} {\node [fill=white] at (\x) {$\x$};}
\node at (5.5,1) {$Hol_1 \cdot \alpha(c^{-1}f^{-1}d)(Hol_2) \cdot \alpha(c^{-1})(Hol_3 \cdot \alpha(h^{-1}k^{-1}i)(Hol_4))$};
\end{tikzpicture}
\end{align*}
whereas if the squares are composed vertically and then horizontally:
\begin{align*}
&\begin{tikzpicture}[baseline={([yshift=-.5ex]current bounding box.center)},vertex/.style={anchor=base,
    circle,fill=black!25,minimum size=18pt,inner sep=2pt}]
\fill (0,4) circle (2pt);
\draw[shorten <= 0.2cm, shorten >= 0.2cm, ->] (0,4) to coordinate[midway]  (a) (2,4);
\draw[shorten <= 0.2cm, shorten >= 0.2cm, ->] (2,4) to coordinate[midway]  (b) (4,4);
\draw[shorten <= 0.2cm, shorten >= 0.2cm, ->] (0,4) to coordinate[midway]  (c) (0,2);
\draw[shorten <= 0.2cm, shorten >= 0.2cm, ->] (2,4) to coordinate[midway]  (d) (2,2);
\draw[shorten <= 0.2cm, shorten >= 0.2cm, ->] (4,4) to coordinate[midway]  (e) (4,2);
\draw[shorten <= 0.2cm, shorten >= 0.2cm, ->] (0,2) to coordinate[midway]  (f) (2,2);
\draw[shorten <= 0.2cm, shorten >= 0.2cm, ->] (2,2) to coordinate[midway]  (g) (4,2);
\draw[shorten <= 0.2cm, shorten >= 0.2cm, ->] (0,2) to coordinate[midway]  (h) (0,0);
\draw[shorten <= 0.2cm, shorten >= 0.2cm, ->] (2,2) to coordinate[midway]  (i) (2,0);
\draw[shorten <= 0.2cm, shorten >= 0.2cm, ->] (4,2) to coordinate[midway]  (j) (4,0);
\draw[shorten <= 0.2cm, shorten >= 0.2cm, ->] (0,0) to coordinate[midway]  (k) (2,0);
\draw[shorten <= 0.2cm, shorten >= 0.2cm, ->] (2,0) to coordinate[midway]  (l) (4,0);
 \foreach \x in {a,b,c,...,l} {\node [fill=white] at (\x) {$\x$};}
 \node at (1,3) {$Hol_1$};
  \node at (3,3) {$Hol_2$};
   \node at (1,1) {$Hol_3$};
    \node at (3,1) {$Hol_4$};
\end{tikzpicture}
= \begin{tikzpicture}[baseline={([yshift=-.5ex]current bounding box.center)},vertex/.style={anchor=base,
    circle,fill=black!25,minimum size=18pt,inner sep=2pt}]
\fill (0,4) circle (2pt);
\draw[shorten <= 0.2cm, shorten >= 0.2cm, ->] (0,4) to coordinate[midway]  (a) (4,4);
\draw[shorten <= 0.2cm, shorten >= 0.2cm, ->] (4,4) to coordinate[midway]  (b) (8,4);
\draw[shorten <= 0.2cm, shorten >= 0.2cm, ->] (0,4) to coordinate[midway]  (hc) (0,0);
\draw[shorten <= 0.2cm, shorten >= 0.2cm, ->] (4,4) to coordinate[midway]  (id) (4,0);
\draw[shorten <= 0.2cm, shorten >= 0.2cm, ->] (8,4) to coordinate[midway]  (je) (8,0);
\draw[shorten <= 0.2cm, shorten >= 0.2cm, ->] (0,0) to coordinate[midway]  (k) (4,0);
\draw[shorten <= 0.2cm, shorten >= 0.2cm, ->] (4,0) to coordinate[midway]  (l) (8,0);
 \foreach \x in {a, b, hc, id, je, k} {\node [fill=white] at (\x) {$\x$};}
 \node at (2,2.5) {$Hol_1\cdot \alpha(c^{-1})(Hol_3) $};
  \node at (6,2.5) {$Hol_2 \cdot \alpha(d^{-1})(Hol_4)$};
\end{tikzpicture}\\
= & \begin{tikzpicture}[ scale=1.2, baseline={([yshift=-.5ex]current bounding box.center)},vertex/.style={anchor=base,
    circle,fill=black!25,minimum size=18pt,inner sep=2pt}]
\fill (0,2) circle (2pt);
\draw[shorten >= 0.2cm, ->] (0,2) to coordinate[midway]  (ba) (11,2) to coordinate[midway] (je) (11,0) to coordinate[midway] (lk) (0,0) to coordinate[midway] (hc) (0,2);
 \foreach \x in {ba, je, lk, hc} {\node [fill=white] at (\x) {$\x$};}
\node at (5.5,1) {$Hol_1\cdot \alpha(c^{-1})(Hol_3) \cdot \alpha((hc)^{-1}k^{-1}(id))(Hol_2 \cdot \alpha(d^{-1})(Hol_4))$};
\end{tikzpicture}
\end{align*}
It is straightforward to check that these two expressions are equal via the crossed module relations and observing the boundaries/targets of the $Hol_{\bullet}$ terms.
\end{proof}
\end{prop}

The following two propositions follow from \cite[Lemmas 2.19 and 2.20]{SWII}, and state how these $Hol$ functions relate to one another.

\begin{prop}[Edge Cube]\label{Edge Cube}
Given a square, $\Sigma$, there is an equation at each intersection $U_{ij}$
\begin{align*}
Hol_{ij}(-,0) = &Hol_{ij}(0,-) \cdot \alpha_{g_{ij}(0,0)^{-1}hol_j(0,-)^{-1}g_{ij}(0,s)hol_i(0,-)}(Hol_i(\Sigma))\\
&\cdot \alpha_{g_{ij}(0,0)^{-1}hol_j(0,-)^{-1}g_{ij}(0,s)}(Hol_{ij}(-,s))\alpha_{g_{ij}(0,0)^{-1}}(Hol_j^{-1}(\Sigma))  \\
&\cdot  \alpha_{g_{ij}(0,0)^{-1}hol_j(-,0)hol_j(t,-)g_{ij}(t,s)hol_i(t,-)}(Hol_{ij}^{-1}(t,-))
\end{align*}
expressed by the following cube-diagram
\begin{equation}
\begin{tikzpicture}[baseline={([yshift=-.5ex]current bounding box.center)},vertex/.style={anchor=base,
    circle,fill=black!25,minimum size=18pt,inner sep=2pt}]
\path[shorten >=0.2cm,shorten <=0.2cm,->]   (3,9)       edge coordinate[midway] (b)node {\contour{white}{$hol_i(-,0)$}}(9,9);
\path[shorten >=0.2cm,shorten <=0.2cm,->]   (3,3)       edge coordinate[midway] (b)node  {\contour{white}{$hol_j(-,0)$}}(9,3);
\path[shorten >=0.2cm,shorten <=0.2cm,->]   (3,9)       edge coordinate[midway] (b)node[rotate=-90] {\contour{white}{$g_{ij}(0,0)$}}(3,3);
\path[shorten >=0.2cm,shorten <=0.2cm,->]   (9,9)       edge coordinate[midway] (b)node[rotate=-90] {\contour{white}{$g_{ij}(t,0)$}}(9,3);	
	\node at (6,6) {$Hol_{ij}(-,0)$};
	\fill (3,9)circle (2pt);
	\fill (9,9)circle (2pt);
	\fill (9,3)circle (2pt);
	\fill (3,3)circle (2pt);
\end{tikzpicture}= 
\begin{tikzpicture}[baseline={([yshift=-.5ex]current bounding box.center)},vertex/.style={anchor=base,
    circle,fill=black!25,minimum size=18pt,inner sep=2pt}]
\path[shorten >=0.2cm,shorten <=0.2cm,->]   (3,9)       edge coordinate[midway] (b)node {\contour{white}{$hol_i(-,0)$}}(9,9);
\path[shorten >=0.2cm,shorten <=0.2cm,->]   (3,3)       edge coordinate[midway] (b)node  {\contour{white}{$hol_j(-,0)$}}(9,3);
\path[shorten >=0.2cm,shorten <=0.2cm,->]   (3,9)       edge coordinate[midway] (b)node[rotate=-90] {\contour{white}{$g_{ij}(0,0)$}}(3,3);
\path[shorten >=0.2cm,shorten <=0.2cm,->]   (9,9)       edge coordinate[midway] (b)node[rotate=-90] {\contour{white}{$g_{ij}(t,0)$}}(9,3);
\path[shorten >=0.2cm,shorten <=0.2cm,->]   (4.5,7.5)       edge coordinate[midway] (b)node {\contour{white}{$hol_i(-,s)$}}(7.5,7.5);
\path[shorten >=0.2cm,shorten <=0.2cm,->]   (4.5,4.5)       edge coordinate[midway] (b)node {\contour{white}{$hol_j(-,s)$}}(7.5,4.5);
\path[shorten >=0.2cm,shorten <=0.2cm,->]   (4.5,7.5)       edge coordinate[midway] (b)node[rotate=-90] {\contour{white}{$g_{ij}(0,s)$}}(4.5,4.5);
\path[shorten >=0.2cm,shorten <=0.2cm,->]   (7.5,7.5)       edge coordinate[midway] (b)node[rotate=-90] {\contour{white}{$g_{ij}(t,s)$}}(7.5,4.5);
\path[shorten >=0.2cm,shorten <=0.2cm,->]   (3,9)       edge coordinate[midway] (b)node[rotate=-45] {\contour{white}{$hol_i(0,-)$}}(4.5,7.5);
\path[shorten >=0.2cm,shorten <=0.2cm,->]   (9,9)       edge coordinate[midway] (b)node[rotate=45] {\contour{white}{$hol_i(t,-)$}}(7.5,7.5);
\path[shorten >=0.2cm,shorten <=0.2cm,->]   (3,3)       edge coordinate[midway] (b)node[rotate=45] {\contour{white}{$hol_j(0,-)$}}(4.5,4.5);
\path[shorten >=0.2cm,shorten <=0.2cm,->]   (9,3)       edge coordinate[midway] (b)node[rotate=-45] {\contour{white}{$hol_j(t,-)$}}(7.5,4.5);	
	\node at (6,6) {$Hol_{ij}(-,s)$};
	\node at (6,8.25) {$Hol_{i}(\Sigma)$};
	\node[rotate=-90] at (8.25,6) {$Hol_{ij}^{-1}(t,-)$};
	\node at (6,3.75) {$Hol_{j}^{-1}(\Sigma)$};
	\node[rotate=-90]  at (3.76,6) {$Hol_{ij}(0,-)$};
	\fill (3,9)circle (2pt);
	\fill (9,9)circle (2pt);
	\fill (9,3)circle (2pt);
	\fill (3,3)circle (2pt);
	\fill (4.5,4.5)circle (2pt);
	\fill (4.5,7.5)circle (2pt);
	\fill (7.5,7.5)circle (2pt);
	\fill (7.5,4.5)circle (2pt);
\end{tikzpicture}
\end{equation}
\end{prop}

\begin{prop}[Vertex Cube]\label{Vertex Cube}
For a path $\gamma$ in $U_{ijkl}$, with $x= \gamma(0)$ and $y=\gamma(1)$, the local data for squares of Section \ref{section square data} satisfies
\begin{align*}Hol_{ijkl}(x) = &Hol_{ij}\cdot \alpha_{g_{ij}^{-1}(x) hol_j g_{ij}(y)hol_i}(Hol_{ik}^{-1}) \cdot \alpha_{g_{ij}^{-1}(x) hol_j^{-1} g_{ij}(y)}(Hol_{ijkl}(y)) \\
&\cdot \alpha_{g_{ij}^{-1}(x)}(Hol_{jl})\cdot \alpha_{g_{ij}^{-1}(x)g_{jl}^{-1}(x)hol_l^{-1}g_{kl}(y)hol_k}(Hol_{kl}^{-1})
\end{align*}
expressed by the following cube-diagram
\begin{equation}
\begin{tikzpicture}[baseline={([yshift=-.5ex]current bounding box.center)},vertex/.style={anchor=base,
    circle,fill=black!25,minimum size=18pt,inner sep=2pt}]
\path[shorten >=0.2cm,shorten <=0.2cm,->]   (3,9)       edge coordinate[midway] (b)node {\contour{white}{$g_{ik}(x)$}}(9,9);
\path[shorten >=0.2cm,shorten <=0.2cm,->]   (3,3)       edge coordinate[midway] (b)node  {\contour{white}{$g_{jl}(x)$}}(9,3);
\path[shorten >=0.2cm,shorten <=0.2cm,->]   (3,9)       edge coordinate[midway] (b)node[rotate=-90] {\contour{white}{$g_{ij}(x)$}}(3,3);
\path[shorten >=0.2cm,shorten <=0.2cm,->]   (9,9)       edge coordinate[midway] (b)node[rotate=-90] {\contour{white}{$g_{kl}(x)$}}(9,3);
	\node at (6,6) {$Hol_{ijkl}(x)$};
	\fill (3,9)circle (2pt);
	\fill (9,9)circle (2pt);
	\fill (9,3)circle (2pt);
	\fill (3,3)circle (2pt);
\end{tikzpicture}= 
\begin{tikzpicture}[baseline={([yshift=-.5ex]current bounding box.center)},vertex/.style={anchor=base,
    circle,fill=black!25,minimum size=18pt,inner sep=2pt}]
\path[shorten >=0.2cm,shorten <=0.2cm,->]   (3,9)       edge coordinate[midway] (b)node {\contour{white}{$g_{ik}(x)$}}(9,9);
\path[shorten >=0.2cm,shorten <=0.2cm,->]   (3,3)       edge coordinate[midway] (b)node  {\contour{white}{$g_{jl}(x)$}}(9,3);
\path[shorten >=0.2cm,shorten <=0.2cm,->]   (3,9)       edge coordinate[midway] (b)node[rotate=-90] {\contour{white}{$g_{ij}(x)$}}(3,3);
\path[shorten >=0.2cm,shorten <=0.2cm,->]   (9,9)       edge coordinate[midway] (b)node[rotate=-90] {\contour{white}{$g_{kl}(x)$}}(9,3);
\path[shorten >=0.2cm,shorten <=0.2cm,->]   (4.5,7.5)       edge coordinate[midway] (b)node {\contour{white}{$g_{ik}(y)$}}(7.5,7.5);
\path[shorten >=0.2cm,shorten <=0.2cm,->]   (4.5,4.5)       edge coordinate[midway] (b)node {\contour{white}{$g_{jl}(y)$}}(7.5,4.5);
\path[shorten >=0.2cm,shorten <=0.2cm,->]   (4.5,7.5)       edge coordinate[midway] (b)node[rotate=-90] {\contour{white}{$g_{ij}(y)$}}(4.5,4.5);
\path[shorten >=0.2cm,shorten <=0.2cm,->]   (7.5,7.5)       edge coordinate[midway] (b)node[rotate=-90] {\contour{white}{$g_{kl}(y)$}}(7.5,4.5);
\path[shorten >=0.2cm,shorten <=0.2cm,->]   (3,9)       edge coordinate[midway] (b)node[rotate=-45] {\contour{white}{$hol_i$}}(4.5,7.5);
\path[shorten >=0.2cm,shorten <=0.2cm,->]   (9,9)       edge coordinate[midway] (b)node[rotate=45] {\contour{white}{$hol_k$}}(7.5,7.5);
\path[shorten >=0.2cm,shorten <=0.2cm,->]   (3,3)       edge coordinate[midway] (b)node[rotate=45] {\contour{white}{$hol_i$}}(4.5,4.5);
\path[shorten >=0.2cm,shorten <=0.2cm,->]   (9,3)       edge coordinate[midway] (b)node[rotate=-45] {\contour{white}{$hol_l$}}(7.5,4.5);	
	\node at (6,6) {$Hol_{ijkl}(y)$};
	\node at (6,8.25) {$Hol_{ik}^{-1}$};
	\node[rotate=90] at (8.25,6) {$Hol_{kl}^{-1}$};
	\node at (6,3.75) {$Hol_{jl}$};
	\node[rotate=90]  at (3.76,6) {$Hol_{ij}$};
	\fill (3,9)circle (2pt);
	\fill (9,9)circle (2pt);
	\fill (9,3)circle (2pt);
	\fill (3,3)circle (2pt);
	\fill (4.5,4.5)circle (2pt);
	\fill (4.5,7.5)circle (2pt);
	\fill (7.5,7.5)circle (2pt);
	\fill (7.5,4.5)circle (2pt);
\end{tikzpicture}
\end{equation}
\end{prop}

\section{Glueing Together Local 2-Holonomy}\label{section glue Hols}
In this section, a definition is provided for 2-holonomy of a square mapped into $M$, $\Sigma \in \mathcal{N}\subset M^{Sq}$ (Definition \ref{def of N}), which lands in multiple open sets $U_i \subset M$.  
Recall that in Section \ref{sec diffeo}, following the notation in \cite{TWZ}, an open cover of $M^{Sq}$ is given by sets $$\mathcal{N}_I= \{\Sigma \in M^{Sq} \vert \text{ for each } (p,q) \in I, \Sigma(Sq_{(p,q)}) \subset U_{i_{(p,q)}} \}.$$
While the notation and inspiration for this particular paper is attributed to \cite{TWZ}, the broader scope of the ideas in this section (in particular the non-abelian approach to glueing squares) can be traced back to \cite{BrHi} while the details and diagrams of \cite{MP2} (eg Fig 3) were not only inspirational to this paper, but serve as a reference for this cubical approach to ``patching together local holonomies''.

\subsection{Comments on the globalness of 2-Holonomy}
Before providing a definition for a ``global'' 2-holonomy in the next section, it is important to address the extent to which it is global.  It is helpful to first review what this means for the 1-holonomy given by a principal $G$-bundle with connection, given by local data subordinate to an open cover, $\{U_i\}$.  The first thing one obtains from this local data is a functor of groupoids, $\mathcal{P}(U_i) \xrightarrow{hol_i} BG$, from the path 1-groupoid referenced above to the groupoid with a single object and morphism space equal to the Lie group, $G$.  The fact that this is a functor implies in particular, that the holonomy of a constant path is the identity, i.e. $hol_i( \text{constant path} ) = id_G$ as well as that holonomy respects composition of paths, i.e. $hol_i (\gamma_1 \circ \gamma_2) = hol_i (\gamma_1) \cdot hol_i( \gamma_2)$.  The space of paths on the base manifold, up to thin homotopy again, has open sets of the form $\mathcal{N}_{\{i_1, i_2, \ldots, i_n\}}$, where similar to definition \ref{def of N}, the paths in this open set are decomposed into pieces which land in the open sets $U_{i_1}, \ldots, U_{i_n}$.  By applying the $g_{ij}$ transition maps at the vertices of these decomposed paths, and multiplying the $hol_{i_k}$ maps to the pieces of path, we obtain maps, $hol_{\mathcal{N}} : \mathcal{N} \to G$.  Next, restricting to loops, we can consider the sub-groupoid of loops $\mathcal{L}(M)$ over the base manifold, $M$, with similar open sets, and see that these maps $hol_{\mathcal{N}}$ differ between two open sets $\mathcal{N}$ and $\mathcal{N}'$ by conjugation of $g_{i_1 i'_1}$ applied to the based point of the loop in question.  Finally, we can therefore say there is a well-defined map of groupoids, $hol: \mathcal{L}(M) \to \mathcal{B}$, where $\mathcal{B}$ is a bundle of groupoids, with trivializations of the form $\mathcal{N} \times BG$.  It is in such a context that one can say that the properties of 1-holonomy can provide a \emph{global map}.

Note however that for 2-holonomy, we may not have a transitive relation when we try to build a bundle of 2-groupoids (since our $g_{ij}$ do not statisfy the cocycle condition).  In future work based on a current joint project with Micah Miller, Thomas Tradler, and Mahmoud Zeinalian, we hope to formally establish the previous paragraph and then setup and prove the following claim: 
\begin{quote}
\emph{The 2-holonomy maps, $Hol_{\mathcal{N}}$, described in definition \ref{def of Hol} glue together to construct a map of 2-groupoids}
\[ Hol: \mathcal{S}_2(M) \to \tilde{\mathcal{G}} \]
\emph{from the sphere 2-groupoid to a gerbe of $2$-groupoids over $M$ with structure 2-groupoid given by the crossed module, $\mathcal{G}$.}
\end{quote}
While the details required to carefully setup and prove all of the statements implied in this subsection are not all contained in this paper, the following subsection goes a long way to collect many of the propositions which would be essential to attaining the global 2-holonomy map described above.

\subsection{Semi-Global 2-Holonomy}
For the remainder of this section, fix a good open cover $\mathcal{U} = \{U_i\}$ of $M$ and an open set $\mathcal{N} \subset M^{Sq}$ as in Definition \ref{def of N}. 
\begin{definition}\label{def of Hol}
Given local transport data for squares, $\{Hol_i\}$, $\{Hol_{ij}\}$, $\{Hol_{ijkl}\}$, define $Hol^{\mathcal{N}}: \mathcal{N} \to H$ on a square $\Sigma \in \mathcal{N}\subset M^{Sq}$ by assembling the local data on the \emph{grid} (Definition \ref{def of grid}) as shown:
\begin{equation}\label{general Hol grid}
 \resizebox{8cm}{!}{ \begin{tikzpicture}
     \foreach \t in {0,2,3.5,5.5,7,9,10,12}
    \foreach \s in {0,2,3,5,6.5,8.5,10,12}
    {
        	\fill (\t,\s)circle (2pt);
    }
            \foreach \t in {1, 2.75, 4.5, 6.25, 8,11}
        {    	\node at (\t,2.5) {$\vdots$};}  
        
             \foreach \s in {1, 4, 5.75, 7.5, 9.25, 11}
        {    	\node at (9.5,\s) {$\dots$};} 
        
     \foreach \t in {0,3.5,7,10}
    \foreach \s in {0,2,3,5,6.5,8.5,10,12}
    {
        \path[shorten >=0.3cm,shorten <=0.3cm,->]  (\t,\s)       edge (\t+2,\s);
    
    }
    
      \foreach \t in {2,5.5}
    \foreach \s in {0,2,3,5,6.5,8.5,10,12}
    {
        \path[shorten >=0.3cm,shorten <=0.3cm,->]  (\t,\s)       edge (\t+1.5,\s);
    
    }

          \foreach \t in {0,2,3.5,5.5,7,9,10,12}
    \foreach \s in {2,5,8.5,12}
    {
        \path[shorten >=0.3cm,shorten <=0.3cm,->]  (\t,\s)       edge (\t,\s-2);
    
    }
              \foreach \t in {0,2,3.5,5.5,7,9,10,12}
    \foreach \s in {6.5,10}
    {
        \path[shorten >=0.3cm,shorten <=0.3cm,->]  (\t,\s)       edge (\t,\s-1.5);
    
    }

        \node at(9.5, 2.5) {$\ddots$};
     \node  at (1, 11) {$Hol_{a}$};
          \node  at (4.5, 11) {$Hol_{b}$};
     \node  at (8, 11) {$Hol_{c}$};
     \node  at (11, 11) {$Hol_{d}$};
     
        \node  at (2.75, 11) {$Hol_{ab}^{-1}$};
          \node  at (6.25, 11) {$Hol_{bc}^{-1}$};
      
      \node  at (1, 9.25) {$Hol_{ae}$};
          \node  at (4.5, 9.25) {$Hol_{bf}$};
     \node  at (8, 9.25) {$Hol_{cg}$};
     \node  at (11, 9.25) {$Hol_{dh}$};
     
        \node  at (2.75, 9.25) {$Hol_{aebf}$};
          \node  at (6.25, 9.25) {$Hol_{bfcg}$};

     \node  at (1, 7.5) {$Hol_{e}$};
          \node  at (4.5, 7.5) {$Hol_{f}$};
     \node  at (8, 7.5) {$Hol_{g}$};
     \node  at (11, 7.5) {$Hol_{h}$};
     
             \node  at (2.75, 7.5) {$Hol_{ef}^{-1}$};
          \node  at (6.25, 7.5) {$Hol_{fg}^{-1}$};

        \node  at (1, 5.75) {$Hol_{ei}$};
          \node  at (4.5, 5.75) {$Hol_{fj}$};
     \node  at (8, 5.75) {$Hol_{gk}$};
     \node  at (11, 5.75) {$Hol_{hl}$};
     
        \node  at (2.75, 5.75) {$Hol_{eifj}$};
          \node  at (6.25, 5.75) {$Hol_{fjgk}$};
          
     \node  at (1, 4) {$Hol_{i}$};
          \node  at (4.5, 4) {$Hol_{j}$};
     \node  at (8, 4) {$Hol_{k}$};
     \node  at (11, 4) {$Hol_{l}$};      
          
   \node  at (2.75, 4) {$Hol_{ij}^{-1}$};
   \node  at (6.25, 4) {$Hol_{jk}^{-1}$};
             
     \node  at (1, 1) {$Hol_{m}$};
          \node  at (4.5, 1) {$Hol_{n}$};
     \node  at (8, 1) {$Hol_{o}$};
     \node  at (11, 1) {$Hol_{p}$};      
          
   \node  at (2.75, 1) {$Hol^{-1}_{mn}$};
   \node  at (6.25, 1) {$Hol^{-1}_{no}$};

  \end{tikzpicture}}
  \end{equation}
Using the multiplication conventions (diagrams \eqref{horizontal glue} and \eqref{vertical glue}) for squares, glue first vertically and then horizontally to obtain the following expression
 \begin{align}\label{eq global Hol}
 Hol^{\mathcal{N}}:= &Hol_a \cdot \overline{Hol_{ae}} \cdot \overline{Hol_e} \cdot \ldots \cdot \overline{Hol_m}\\
 \cdot &\overline{Hol^{-1}_{ab}} \cdot \overline{Hol_{aebf}} \cdot \overline{Hol^{-1}_{ef}} \cdot \ldots \cdot \overline{Hol^{-1}_{mn}} \nonumber \\
 \vdots & \nonumber \\
 \cdot & \overline{Hol_d} \cdot \overline{Hol_{dh}} \cdot \overline{Hol_h} \cdot \cdot \ldots \cdot \overline{Hol_p}\nonumber
 \end{align}
 where $Hol_i$ is evaluated on the \emph{face} $\Sigma_i$, $Hol_{ij}$ is evaluated on the \emph{edge} $\gamma_{ij}$, and $Hol_{ijkl}$ is evaluated on the \emph{vertex} $x_{ijkl}$ as described in Definition \ref{def of grid}.  The overline decoration on each, ``$\overline{X}$'', is described in Convention \ref{conv overline} below.  Note that each face, $\Sigma_a$, and edge, $\gamma_{ab}$, technically corresponds to (products of) subintervals of $I$.  However, note that reparametrization does not change the output since these functions are based on the thin-homotopy-invariant functions of \cite{SWII}.
\end{definition} 

\begin{convention}\label{conv overline}
In general, the overline decoration is shorthand for an $\alpha$-action, where $\overline{X} = \alpha_g(X)$ for the appropriate element, $g \in G$, which we now explain.  The rules for glueing horizontally and vertically (diagrams \eqref{horizontal glue} and \eqref{vertical glue}) determine the element, $g$, once the order of the squares is chosen.  When the order is not explicitly stated, the convention of first glueing vertically and then horizontally is assumed.  In turn, the action on $X$ is given by the following 1-holonomy: start at the upper left corner of the grid, then move down the far left edge of the grid until you reach the bottom, then move right until you are at the left side of the column of $X$, then move up until you end at the upper left corner of the $X$-square as seen below: 
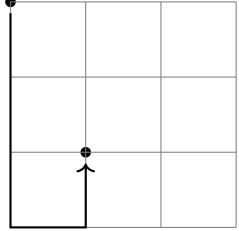
\begin{equation}\label{pic of path}  \begin{tikzpicture}
	\fill (1,1) circle (2pt);
  	\fill (0,3) circle (2pt);
     \foreach \t in {0,1, 2, 3}
     {
        \path[thin, gray] (\t, 3) edge (\t, 0);
    }
    \foreach \s in {0,1, 2, 3}
    {
        \path[thin, gray] (0, \s) edge (3, \s);
    }

        \draw[thick, shorten >=0.15cm,shorten <=0.15cm,->]  (0,3) -- (0,0) -- (1,0) -- (1,1);
    
    \end{tikzpicture} 
\end{equation}
For example, in \eqref{horizontal glue}, we would write $\overline{Hol_2} := \alpha_{d^{-1} c^{-1} b} (Hol_2)$, while in \eqref{vertical glue} we would write $\overline{Hol_2} := \alpha_{a^{-1}}(Hol_2)$.
\end{convention}

\subsection{Properties of $Hol^{\mathcal{N}}$ for Squares}\label{sec props of Hol}
This section will provide a summary of the properties one would expect $Hol^{\mathcal{N}}$ to have.  Note that the propositions here are mostly modifications to the works of Schreiber and Waldorf along with the works of Martins and Picken.  In places where the proof might be original to the context of this paper, some additional comments are supplied.

\begin{prop}\label{prop Hol invar thin homotopy}
$Hol^{\mathcal{N}}$ is invariant under thin homotopy.
\begin{proof}
See \cite[Lemma 2.16]{SWII}.
\end{proof}
\end{prop}

\begin{prop}\label{prop t of Hol and composition}
The target of $Hol^{\mathcal{N}}$ is equal to the one-holonomy along the boundary of the grid.  Explicitly, assuming the grid \eqref{general Hol grid}, is
$$t(Hol^{\mathcal{N}}) = hol_a^{-1} \cdot g_{ab}^{-1} \cdot hol_b^{-1} \cdot \ldots \cdot hol_d^{-1} \cdot g_{dh}^{-1} \cdot \ldots hol_o \cdot g_{no}\cdot hol_n \cdot \ldots g_{ei} \cdot hol_e \cdot g_{ae} \cdot hol_a.$$
Furthermore, $Hol^{\mathcal{N}}$ respects composition of squares.
\end{prop}

\begin{prop}\label{prop Hol invar subdivision}
$Hol^{\mathcal{N}}$ is invariant under subdivision.
\begin{proof}
Subdivision induces a new grid, on which the added transition data will be shown to be the identity.  Consider the following vignette of what subdivision might produce near an $ijkl$ vertex:
\begin{equation}
 \resizebox{6cm}{!}{\begin{tikzpicture}[baseline={([yshift=-.5ex]current bounding box.center)},vertex/.style={anchor=base,
     circle,fill=black!25,minimum size=18pt,inner sep=2pt}]   
      \foreach \t in {0,4,6,10} {  \path  (\t, 10)       edge (\t,0);}
 \foreach \s in {0,4,6,10}{  \path  (10, \s)       edge (0,\s);}
\node  at (2, 8) {$Hol_i$};
\node  at (2, 5) {$Hol_{ij}$};
\node  at (2, 2) {$Hol_{j}$};
\node  at (5, 8) {$Hol_{ik}^{-1}$};
\node  at (5, 5) {$Hol_{ijkl}$};
\node  at (5, 2) {$Hol_{jl}^{-1}$};
\node  at (8, 8) {$Hol_k$};
\node  at (8, 5) {$Hol_{kl}$};
\node  at (8, 2) {$Hol_l$};  
     \end{tikzpicture}}
     = 
 \resizebox{6cm}{!}{\begin{tikzpicture}[baseline={([yshift=-.5ex]current bounding box.center)},vertex/.style={anchor=base,
     circle,fill=black!25,minimum size=18pt,inner sep=2pt}]   
 \foreach \t in {0,4,6,10} {  \path[thick]  (\t, 10)       edge (\t,0);}
 \foreach \s in {0,4,6,10}{  \path[thick]  (10, \s)       edge (0,\s);}
  \foreach \t in {1.3,2.7} {  \path[thin, gray]  (\t, 10)       edge (\t,0);}
 \foreach \s in {7.5,8.5}{  \path[thin, gray]  (10, \s)       edge (0,\s);}
\node  at (0.75, 9.25) {$Hol_i$};
\node  at (0.75, 8) {$Hol_i$};
\node  at (0.75, 6.75) {$Hol_i$};
\node  at (0.75, 5) {$Hol_{ij}$};
\node  at (0.75, 2) {$Hol_{j}$};

\node  at (2, 9.25) {$Hol^{-1}_{ii}$};
\node  at (2, 8) {$Hol_{iiii}$};
\node  at (2, 6.75) {$Hol^{-1}_{ii}$};
\node  at (2, 5) {$Hol_{ijij}$};
\node  at (2, 2) {$Hol^{-1}_{jj}$};

\node  at (3.25, 9.25) {$Hol_i$};
\node  at (3.25, 8) {$Hol_i$};
\node  at (3.25, 6.75) {$Hol_i$};
\node  at (3.25, 5) {$Hol_{ij}$};
\node  at (3.25, 2) {$Hol_{j}$};

\node  at (5, 9.25) {$Hol_{ik}^{-1}$};
\node  at (5, 8) {$Hol_{iikk}$};
\node  at (5, 6.75) {$Hol_{ik}^{-1}$};
\node  at (5, 5) {$Hol_{ijkl}$};
\node  at (5, 2) {$Hol_{jl}^{-1}$};

\node  at (8, 9.25) {$Hol_k$};
\node  at (8, 8) {$Hol_{kk}$};
\node  at (8, 6.75) {$Hol_{k}$};
\node  at (8, 5) {$Hol_{kl}$};
\node  at (8, 2) {$Hol_{l}$};   
     \end{tikzpicture}}
\end{equation}
where equality comes from $Hol_{ii}$, $Hol_{iiii}$, and $Hol_{ijij}$ all being the identity in $H$.
\end{proof}
\end{prop}

\begin{prop}\label{prop Hol transforms}
$Hol^{\mathcal{N}}$ transforms across open sets $\mathcal{N}_I$ and $\mathcal{N}_{I'}$ by $\alpha_{g_{ii'}}$ at the base point and by $Hol_{jj'}$ and $Hol_{iji'j'}$ along the boundary.  In particular:
$$Hol^{\mathcal{N}_{I'}} = \alpha_{g_{ii'}}(Hol^{\mathcal{N}_{I}}) \cdot \prod\limits_{\partial \Sigma} \overline{Hol_{jj'}} \cdot \prod\limits_{\partial \Sigma} \overline{Hol_{jkj'k'}}$$
\begin{proof}
Let $\Sigma \in \mathcal{N}_{I_0} \cap \mathcal{N}_{I_0'}$.   By  Proposition \ref{prop Hol invar subdivision}, there exists a subdivision using open sets $\mathcal{N}_{I}$ and $\mathcal{N}_{I'}$ which both use a grid of size $n$ ($t$-direction) by $m$ ($s$-direction) such that 
$$Hol^{\mathcal{N}_{I_0}} = Hol^{\mathcal{N}_{I}}  \quad \text{and} \quad Hol^{\mathcal{N}_{I'_0}} = Hol^{\mathcal{N}_{I'}}.$$
Consider the grid of $Hol^{\mathcal{N}_{I'}}$ analogous to that from Definition \ref{def of Hol} but replace each face $Hol_{i'}$ with a cube (Proposition \ref{Edge Cube}), replace each horizontal edge $i'j'$ with a cube (Proposition \ref{Vertex Cube}), and similarly each vertical edge $i'k'$ with a cube.  Although it has not been laid out in a previous proposition, one can divide an $ijkli'j'k'l'$ cube into $i_1 i_2 i_3 i_4$-tetrahedra to build a cube: 
  \begin{equation}
 \resizebox{8cm}{!}{\begin{tikzpicture}[scale=0.8, baseline={([yshift=-.5ex]current bounding box.center)},vertex/.style={anchor=base,
    circle,fill=black!25,minimum size=18pt,inner sep=2pt}]
\path[shorten >=0.2cm,shorten <=0.2cm,->]   (3,9)       edge coordinate[midway] (b)node {}(9,9);
\path[shorten >=0.2cm,shorten <=0.2cm,->]   (3,3)       edge coordinate[midway] (b)node  {}(9,3);
\path[shorten >=0.2cm,shorten <=0.2cm,->]   (3,9)       edge coordinate[midway] (b)node[rotate=-90] {}(3,3);
\path[shorten >=0.2cm,shorten <=0.2cm,->]   (9,9)       edge coordinate[midway] (b)node[rotate=-90] {}(9,3);	
	\node at (6,6) {$Hol_{i'j'k'l'}$};
	\fill (3,9)circle (2pt);
	\node[above left] at (3,9) {$(0,0)$};
	\fill (9,9)circle (2pt);
	\node[above right] at (9,9) {$(t,0)$};
	\fill (9,3)circle (2pt);
	\node[below right] at (9,3) {$(s,t)$};	
	\fill (3,3)circle (2pt);
	\node[below left] at (3,3) {$(0,s)$};
	
\end{tikzpicture} = 
\begin{tikzpicture}[scale=0.8, baseline={([yshift=-.5ex]current bounding box.center)},vertex/.style={anchor=base,
    circle,fill=black!25,minimum size=18pt,inner sep=2pt}]
\path[shorten >=0.2cm,shorten <=0.2cm,->]   (3,9)       edge coordinate[midway] (b)node {}(9,9);
\path[shorten >=0.2cm,shorten <=0.2cm,->]   (3,3)       edge coordinate[midway] (b)node  {}(9,3);
\path[shorten >=0.2cm,shorten <=0.2cm,->]   (3,9)       edge coordinate[midway] (b)node[rotate=-90] {}(3,3);
\path[shorten >=0.2cm,shorten <=0.2cm,->]   (9,9)       edge coordinate[midway] (b)node[rotate=-90] {}(9,3);
\path[shorten >=0.2cm,shorten <=0.2cm,->]   (4.5,7.5)       edge coordinate[midway] (b)node {}(7.5,7.5);
\path[shorten >=0.2cm,shorten <=0.2cm,->]   (4.5,4.5)       edge coordinate[midway] (b)node {}(7.5,4.5);
\path[shorten >=0.2cm,shorten <=0.2cm,->]   (4.5,7.5)       edge coordinate[midway] (b)node[rotate=-90] {}(4.5,4.5);
\path[shorten >=0.2cm,shorten <=0.2cm,->]   (7.5,7.5)       edge coordinate[midway] (b)node[rotate=-90] {}(7.5,4.5);
\path[shorten >=0.2cm,shorten <=0.2cm,<-]   (3,9)       edge coordinate[midway] (b)node[rotate=-45] {}(4.5,7.5);
\path[shorten >=0.2cm,shorten <=0.2cm,<-]   (9,9)       edge coordinate[midway] (b)node[rotate=45] {}(7.5,7.5);
\path[shorten >=0.2cm,shorten <=0.2cm,<-]   (3,3)       edge coordinate[midway] (b)node[rotate=45] {}(4.5,4.5);
\path[shorten >=0.2cm,shorten <=0.2cm,<-]   (9,3)       edge coordinate[midway] (b)node[rotate=-45] {}(7.5,4.5);	
	\node at (6,6) {$Hol_{ijkl}$};
	\node at (6,8.25) {$Hol_{iki'k'}$};
	\node[rotate=-90] at (8.25,6) {$Hol_{klk'l'}$};
	\node at (6,3.75) {$Hol_{jlj'l'}^{-1}$};
	\node[rotate=-90]  at (3.76,6) {$Hol_{iji'j'}^{-1}$};
	\fill (3,9)circle (2pt);
	\fill (9,9)circle (2pt);
	\fill (9,3)circle (2pt);
	\fill (3,3)circle (2pt);
	\fill (4.5,4.5)circle (2pt);
	\fill (4.5,7.5)circle (2pt);
	\fill (7.5,7.5)circle (2pt);
	\fill (7.5,4.5)circle (2pt);
\end{tikzpicture}}
\end{equation}
Replacing each function in the grid of $Hol^{\mathcal{N}_{I'}}$ with their associated cube from above, all of the interior transition data cancels, leaving only transition data at the boundary.  Note that the $\alpha_{g_{ii'}}$ from the statement comes from changing the basepoint from $U_i$ to $U_{i'}$ in the upper left corner of the square.  This proves the statement of the proposition.
\end{proof}
\end{prop}

\section{$d(Hol)$ for Squares}\label{sec dHol}
The focus of this section, is the proof of Theorem \ref{dhol nice}, the main result of this paper, which says that the deRham differential applied to global 2-holonomy amounts to replacing one $B_{i}$ in any summand with one $H_i$; in addition to some terms associated to the boundary of $\Sigma$.  

\subsection{A warmup example for the Main Theorem}
Before stating the main theorem of this paper formally in the next section, a warmup to the notation and the main idea is provided here.  Consider a good open cover, $\mathcal{U} := \{ U_1, U_2, U_3, U_4\}$, whose \v{C}ech nerve contains the cells given by the diagram illustrated below:
\begin{equation}
\begin{tikzpicture}[baseline={([yshift=-.5ex]current bounding box.center)},vertex/.style={anchor=base,
    circle,fill=black!25,minimum size=18pt,inner sep=2pt}, scale=0.6]
\fill (0,4) circle (2pt); \node [above left] at (0,4) {$U_1$};
\fill (4,4) circle (2pt); \node [above right] at (4,4) {$U_3$};
\fill (4,0) circle (2pt); \node [below left] at (4,0) {$U_4$};
\fill (0,0) circle (2pt); \node [below left] at (0,0) {$U_{2}$};
\draw[shorten <= 0.2cm, shorten >= 0.2cm, -] (0,4) to coordinate[midway]  (a) (0,0);
\draw[shorten <= 0.2cm, shorten >= 0.2cm, -] (0,0) to coordinate[midway]  (b) (4,0);
\draw[shorten <= 0.2cm, shorten >= 0.2cm, -] (4,0) to coordinate[midway]  (c) (4,4);
\draw[shorten <= 0.2cm, shorten >= 0.2cm, -] (0,4) to coordinate[midway]  (d) (4,4);
\draw[shorten <= 0.2cm, shorten >= 0.2cm, -] (4,4) to coordinate[midway]  (e) (0,0);
 \node [fill=white] at (a) {$U_{12}$};
 \node[fill=white]  at (b) {$U_{24}$};
 \node[fill=white]  at (c) {$U_{34}$};
 \node[fill=white]  at (d) {$U_{13}$};
 \node at (1,3) {$U_{123}$};
\node at (3,1) {$U_{234}$};
\end{tikzpicture}
\end{equation}
Given a non-abelian gerbe with connection (Definition \ref{def gerbe}), $\left( \{ g_{ij} \}, \{f_{ijk} \}, \{A_i\}, \{a_{ij}, \{B_i\}\} \right)$, and a one-parameter family of squares, $\Sigma_r \in \mathcal{N}$ (Definition \ref{def of grid}), where for our open sets we have $\mathcal{N} = \mathcal{N}_I$, so that $I = \{ 1,2 \} \times \{1,2\}$ with $U_{(1,1)}:= U_1$, $U_{(1,2)}:= U_2$, $U_{(2,1)}:= U_3$, and $U_{(2,2)}:= U_4$.

By definition \ref{def of Hol}, we have for each $r \in \mathbb{R}$, an element, $Hol^{\mathcal{N}}\left(\Sigma_r\right) \in H$, using the appropriately modified version of that definition's diagram below:
\begin{equation}
\begin{tikzpicture}[baseline={([yshift=-.5ex]current bounding box.center)},vertex/.style={anchor=base,
    circle,fill=black!25,minimum size=18pt,inner sep=2pt}, scale=1.3]
     \foreach \t in {0,2,3,5}
    \foreach \s in {0,2,3,5}
    {
        	\fill (\t,\s)circle (2pt);
    }        
     \foreach \t in {0,3}
    \foreach \s in {0,2,3,5}
    {
        \path[shorten >=0.2cm,shorten <=0.2cm,->]  (\t,\s)       edge (\t+2,\s);    
    }
    
    \foreach \s in  {0,2,3,5}
    {
        \path[shorten >=0.2cm,shorten <=0.2cm,->]  (2,\s)       edge (3,\s);   
    }
     
          \foreach \t in  {0,2,3,5}
    \foreach \s in {2,5}
    {
        \path[shorten >=0.2cm,shorten <=0.2cm,->]  (\t,\s)       edge (\t,\s-2);    
    }
              \foreach \t in {0,2,3,5}
    {
        \path[shorten >=0.2cm,shorten <=0.2cm,->]  (\t,3)       edge (\t,2);    
    }
\node  at (1, 4) {$Hol_1$}; \node  at (1, 1) {$Hol_2$};  \node  at (4, 4) {$Hol_3$};  \node  at (4, 1) {$Hol_4$};
\node  at (1, 2.5) {$Hol_{12}$};  \node  at (2.5, 4) {$Hol_{13}$};  \node  at (4, 2.5) {$Hol_{34}$};  \node  at (2.5, 1) {$Hol_{24}$}; 
\node at (2.5, 2.5) {$Hol_{1234}$};
\end{tikzpicture}
\end{equation}
The goal of Theorem \ref{dhol nice} is to offer a formula for $\restr{\frac{\partial}{\partial r} \left( Hol^{\mathcal{N}} \left( \Sigma_r \right) \right)}{r=0}$.  Recall that the left action for our crossed module of Lie Groups, $\alpha: G \times H \to H$ induces a map $\alpha_h: G \to H$, for each choice of $h \in H$.  In turn, we obtain from the pushforward a map $(\alpha_h)_*: T_g G \to T_{\alpha_g(h)} H$ and thus a map $(\alpha_h)_*: \mathfrak{g} \to T_h H$.  Similarly, $\alpha$ induces for each fixed $g \in G$, a map $\alpha_g: H \to H$ whose pushforward gives a map $(\alpha_g)_*: \mathfrak{h} \to \mathfrak{h}$.  For a fixed $h \in H$ we also have the pushforward of, for example, multiplicatoin on the left, $(L_h)_*: T_{h'} H \to T_{h\cdot h'} H$ and thus $(L_h)_*: \mathfrak{h} \to T_h H$.

Finally, we are ready to examine the content of Theorem \ref{dhol nice} in our example with four open sets.  As stated below, the abbreviated version of the formula for $d(Hol)$ is 
\begin{equation*}
d(Hol)=   - (\alpha_{Hol})_*(A_{i_{(1,1)}}) +Hol\cdot \int_{Sq} H+  Hol \cdot \left(\int_{\partial Sq} B  +\sum\limits_{\partial Sq} a \right) \end{equation*}
The term, $ - (\alpha_{Hol})_*(A_{i_{(1,1)}}) =  - (\alpha_{Hol})_*(A_{1})$ is the pushforward $-(\alpha_{Hol})_*$ applied to $A_{1}\left( \restr{\frac{\partial}{\partial r} \Sigma_r (0,0)}{r=0} \right).$  The term, $\int_{Sq} H$ is a sum of four terms, one for each $3$-form, $H_i \in \Omega^3(U_i, \mathfrak{h})$:
\[ \int_{Sq} H = \int\limits_{Sq_{(1,1)}} (\alpha_{path_{(1,1)}})_* H_1 +  \int\limits_{Sq_{(1,2)}} (\alpha_{path_{(1,2)}})_* H_2 +  \int\limits_{Sq_{(2,2)}} (\alpha_{path_{(2,2)}})_* H_3 +  \int\limits_{Sq_{(2,2)}} (\alpha_{path_{(2,2)}})_* H_4.\]
To be precise, and using the \emph{north, south, east, east} notation of  Definition \ref{def of grid}, we can use the $(2,2)$-square as an example to note that
\begin{align*}
&\int\limits_{Sq_{(2,2)}} \left( \alpha_{path_{(2,2)}}\right)_* H_4\\
 = &\left( \alpha_{hol_1^{-1}(\lambda_1^W) \cdot g^{-1}_{12}(x_{(1,2)}^0) \cdot hol_2^{-1}(\lambda_2^W) \cdot  hol_2^{-1}(\lambda_2^S) \cdot  g^{-1}_{24}(x_{(2,4)}^1) \cdot hol_4(\lambda_{24})} \right)_* \left(\int\limits_{Sq} (\alpha_{path_{(s,t)}})_* (H_4) \right),
\end{align*} 
where $path_{(-,-)} : [0,1]^2 \to G$ is locally defined so that $path_{(s,t)}$ is the $1$-holonomy vertically ``down'' to level $s$ and then horizontally ``right'' over to time, $t$.  Note that one can show the choice of path (as long as it is smooth) does not matter since $t(H_4) = 0$.  Since $\int\limits_{Sq_(2,2)} (\alpha_{path})_* H_4$ is valued in $\mathfrak{h}$, then the pushforward of left-multiplication by $Hol$  means that the term, $Hol \cdot \int_{Sq} H$, takes values in $T_{Hol} H$.
Now that the notation has been explored a bit further, the meaning of the boundary terms becomes clearer, of which there are \emph{twelve}:
\begin{equation}
\begin{tikzpicture}[baseline={([yshift=-.5ex]current bounding box.center)},vertex/.style={anchor=base,
    circle,fill=black!25,minimum size=18pt,inner sep=2pt}, scale=3]
     \foreach \t in {0,1,2,3}
    \foreach \s in {0,1,2,3}
    {
        	\fill (\t,\s)circle (1pt);
    }        
     \foreach \t in {0,2}
    \foreach \s in {0,1,2,3}
    {
        \path[shorten >=0.2cm,shorten <=0.2cm,->]  (\t,\s)       edge (\t+1,\s);    
    }
    
    \foreach \s in  {0,1,2,3}
    {
        \path[shorten >=0.2cm,shorten <=0.2cm,->]  (1,\s)       edge (2,\s);   
    }
     
          \foreach \t in  {0,1,2,3}
    \foreach \s in {2,3}
    {
        \path[shorten >=0.2cm,shorten <=0.2cm,->]  (\t,\s)       edge (\t,\s-1);    
    }
              \foreach \t in {0,1,2,3}
    {
        \path[shorten >=0.2cm,shorten <=0.2cm,->]  (\t,1)       edge (\t,0);    
    }
\node  at (0.5, 2.5) {\color{gray}{$Hol_1$}}; \node  at (0.5, 0.5) {\color{gray}{$Hol_2$}};  \node  at (2.5, 2.5) {\color{gray}{$Hol_3$}};  \node  at (2.5, 0.5) {\color{gray}{$Hol_4$}};
\node  at (0.5, 1.5) {\color{gray}{$Hol_{12}$}};  \node  at (1.5,2.5) {\color{gray}{$Hol_{13}$}};  \node  at (2.5, 1.5) {\color{gray}{$Hol_{34}$}};  \node  at (1.5, 0.5) {\color{gray}{$Hol_{24}$}}; 
\node at (1.5, 1.5) {\color{gray}{$Hol_{1234}$}};
\node[above] at (0.5,3) {$\int_{\gamma^N_{1}} (\alpha_{path^N_{1}})_*(B_{1})$};
\node[above] at (1.5,3) {$(\alpha_{path^N_{13}})_*(a_{13})$};
\node[above] at (2.5,3) {$\int_{\gamma^N_{3}} (\alpha_{path^N_{3}})_*(B_{3})$};
\node[above, rotate=-90] at (3,2.5) {$\int_{\gamma^E_{3}} (\alpha_{path^E_{3}})_*(B_{3})$};
\node[above, rotate=-90] at (3,1.5) {$(\alpha_{path^E_{34}})_*(a_{34})$};
\node[above, rotate=-90] at (3,0.5) {$\int_{\gamma^E_{4}} (\alpha_{path^E_{4}})_*(B_{4})$};
\node[below] at (2.5,0) {$\int_{\gamma^S_{4}} (\alpha_{path^S_{4}})_*(B_{4})$};
\node[below] at (1.5,0) {$(\alpha_{path^S_{24}})_*(a_{24})$};
\node[below] at (0.5,0) {$\int_{\gamma^S_{2}} (\alpha_{path^S_{2}})_*(B_{2})$};
\node[above, rotate=90] at (0,0.5) {$\int_{\gamma^W_{2}} (\alpha_{path^W_{2}})_*(B_{2})$};
\node[above, rotate=90] at (0,1.5) {$(\alpha_{path^W_{12}})_*(a_{12})$};
\node[above, rotate=90] at (0,2.5) {$\int_{\gamma^W_{1}} (\alpha_{path^W_{1}})_*(B_{1})$};
\end{tikzpicture}
\end{equation}

\subsection{The Main Theorem}
For the remainder of this section, a grid, $\mathcal{N}$ is once again fixed.
\begin{theorem}\label{dhol nice}
For the function $Hol=Hol^\mathcal {N}:\mathcal{N}\to H\subset Mat$ defined\footnote{Where $Mat$ is some matrix algebra.} on the open set $\mathcal{N} \subset M^{Sq}$, we can compute the derivative $d (Hol)\in \Omega^1(M^{Sq},Mat)$ as follows:
\begin{equation}\label{eq dHol nice}
d(Hol)=   - (\alpha_{Hol})_*(A_{i_{(1,1)}}) +Hol\cdot \int_{Sq} H+  Hol \cdot \left(\int_{\partial Sq} B  +\sum\limits_{\partial Sq} a \right) \end{equation}
where $i_{(1,1)}$ is the index of the upper left open set $U_{i_{(1,1)}}$, and 
\begin{align}\int_{Sq} H:= & \sum\limits_{\substack{k =1,\dots, n \\ l =1,\dots, m}} \int_{Sq_{(k,l)}} (\alpha_{path_{k,l}})_*(H_{i_{(k,l)}}). \\
\intertext{While the terms in the parenthesis on the right in \eqref{eq dHol nice} are boundary terms:}
\int_{\partial Sq} B:= & - \sum\limits_{k =1,\dots, n} \int_{\gamma^N_{(k,1)}} (\alpha_{path^N_{(k,1)}})_*(B_{i_{(k,1)}}) -  \sum\limits_{l =1,\dots, m} \int_{\gamma^E_{(n,l)}} (\alpha_{path^E_{(n,l)}})_*(B_{i_{(n,l)}})  \\
 &+ \sum\limits_{k=n,\dots, 1} \int_{\gamma^S_{(k,m)}} (\alpha_{path^S_{(k,m)}})_*(B_{i_{(k,m)}}) + \sum\limits_{l =m,\dots, 1} \int_{\gamma^W_{(1,l)}} (\alpha_{path^W_{(1,l)}})_*(B_{i_{(1,l)}}) \nonumber \\
\sum\limits_{\partial Sq} a :=&-  \sum\limits_{k =1,\dots, n-1}(\alpha_{path^N_{(k,1)(k+1,1)}})_*(a_{i_{(k,1)}i_{(k+1,1)}})\nonumber \\
&- \sum\limits_{l =1,\dots, m-1} (\alpha_{path^E_{(n,l)(n,l+1)}})_*(a_{i_{(n,l)}i_{(n,l+1)}}) \nonumber \\
 &+ \sum\limits_{k=n,\dots, 2}  (\alpha_{path^S_{(k-1,m)(k,m)}})_*(a_{i_{(k-1,m)}i_{(k,m)}}) \nonumber \\
  &+ \sum\limits_{l =m,\dots, 2}  (\alpha_{path^W_{(1,l-1)(1,l)}})_*(a_{i_{(1,l-1)}i_{(1,l)}}) \nonumber
 \end{align} 
 where the expressions $\alpha_{path_{\bullet}}(X)$ use an element $path_{\bullet} \in G$ which is given by the appropriate path through the grid.
\end{theorem}
First, a lemma which can be found, for example, in Theorem 2.30 of \cite{MP2}, and which is the local version of this paper's Theorem \ref{dhol nice}, is recalled:

\begin{lemma}\label{local lemma}
For the local 2-holonomy function $Hol_i : U_i^{Sq} \to H$ as defined in Section \ref{section square data} we have
\begin{align}
d(Hol_i)&= - (\alpha_{Hol_i})_*(A_i(0,0)) + Hol_i \cdot \int\limits_{Sq} \alpha_*(H_i) + Hol_i \cdot \int\limits_{\partial Sq} \alpha_*(B_i)
\end{align}
\end{lemma}

\subsection{The Edge and Vertex Relations}\label{Edge and Vertex Relations}
Here two technical lemmas are proved, which give relations amongst the edges and vertices for when $\partiald{r}(Hol)$ is later computed.  Recall the expression for 2-holonomy, Definition \ref{def of Hol} :
 \begin{align*}
 Hol^{\mathcal{N}}= &Hol_a \cdot \overline{Hol_{ae}} \cdot \overline{Hol_e} \cdot \ldots \cdot \overline{Hol_m}\\
 \cdot &\overline{Hol_{ab}^{-1}} \cdot \overline{Hol_{aebf}} \cdot \overline{Hol_{ef}^{-1}} \cdot \ldots \cdot \overline{Hol_{mn}^{-1}}\\
 \vdots &\\
 \cdot & \overline{Hol_d} \cdot \overline{Hol_{dh}} \cdot \overline{Hol_h} \cdot \cdot \ldots \cdot \overline{Hol_p}.
 \end{align*}
 Applying the Leibniz rule, the derivative of global holonomy will start with the expression,
 \begin{align}\label{Liebniz of Hol}
\restr{\partiald{r}}{r=0}Hol : = &\restr{\partiald{r}}{r=0}\left(\ldots \cdot Hol_i \cdot \overline{Hol_{ij} }\cdot \overline{Hol_j} \cdot \overline{Hol_{ik}^{-1}} \cdot \ldots \cdot \overline{Hol_l} \cdot \ldots \right)\\
=& \ldots + \ldots  \cdot \restr{\partiald{r}}{r=0}Hol_i \cdot \overline{Hol_{ij} }\cdot \overline{Hol_j} \cdot \overline{Hol_{ik}^{-1}} \cdot \ldots \cdot \overline{Hol_l} \cdot \ldots \nonumber \\
&+  \ldots \cdot Hol_i \cdot \restr{\partiald{r}}{r=0} \overline{Hol_{ij} }\cdot \overline{Hol_j} \cdot \overline{Hol_{ik}^{-1}} \cdot \ldots \cdot \overline{Hol_l} \cdot \ldots \nonumber \\
 &\vdots  \nonumber \\
 &+ \ldots \cdot Hol_i \cdot  \overline{Hol_{ij} }\cdot \overline{Hol_j} \cdot \overline{Hol_{ik}^{-1}} \cdot \ldots \cdot \restr{\partiald{r}}{r=0}\overline{Hol_l}\cdot \ldots  \nonumber \\
 &+ \ldots.  \nonumber
\end{align}
To state the first lemma, the following setup is used:\\

Assume for the moment that the open set $\mathcal{N}$ only requires four open sets $U_i, U_j, U_k$, and $U_l$, written in counterclockwise order from the upper left of $\Sigma$.  Now, suppose a one-parameter family of squares, $\Sigma(r)$, is given with $\Sigma(0) = \Sigma$.  In particular, then, given the grid associated to $\mathcal{N}$, there is a path, $\rho$, through the vertex $x_{ijkl}$.  Momentarily writing $\rho(0) = x$ and $\rho(r) = y$, there is a vertex cube associated to $\rho$, given by Proposition \ref{Vertex Cube} which corresponds to the equation:
\begin{align*}Hol_{ijkl}(x) = &Hol_{ij}\cdot \alpha_{g_{ij}^{-1}(x) hol_j g_{ij}(y)hol_i}(Hol_{ik}^{-1}) \cdot \alpha_{g_{ij}^{-1}(x) hol_j^{-1} g_{ij}(y)}(Hol_{ijkl}(y)) \\
&\cdot \alpha_{g_{ij}^{-1}(x)}(Hol_{jl})\cdot \alpha_{g_{ij}^{-1}(x)g_{jl}^{-1}(x)hol_l^{-1}g_{kl}(y)hol_k}(Hol_{kl}^{-1}).
\end{align*}  
Thus $Hol_{ijkl}(x)$ in $Hol$ can be replaced, to obtain the following equality of local transport data:
\begin{align}\label{vertex cube}
\resizebox{6cm}{!}{\begin{tikzpicture}[baseline={([yshift=-.5ex]current bounding box.center)},vertex/.style={anchor=base,
    circle,fill=black!25,minimum size=18pt,inner sep=2pt}]
\draw (9,3) -| (12,0)
	node[pos=0.25] {}
	node[pos=0.75, rotate=-90] {\contour{white}{$hol_{l}(\gamma_l^E)$}}
		-| (9,3) 
		node[pos=.25] {\contour{white}{$hol_l(\gamma_l^S)$}}
		node[ pos=.75, rotate=-90] {};	
	\node at (10.5,1.5) {$Hol_{l}(\Sigma_l(0))$};
\draw (9,9) -| (12,3)
	node[pos=0.25] {}
	node[pos=0.75, rotate=-90] {\contour{white}{$g_{kl}(\gamma^h_{kl}(1))$}}
		-| (9,9) 
		node[pos=.25] {\contour{white}{$hol_l(\gamma^h_{kl})$}}
		node[ pos=.75, rotate=-90] { };	
	\node at (10.5,6) {$Hol_{kl}(\gamma^h_{kl})$};
\draw (9,12) -| (12,9)
	node[pos=0.25] {\contour{white}{$g_{ik}(\gamma_{ik}^v(0))$}}
	node[ pos=0.75,rotate=-90] {\contour{white}{$hol_k(\gamma_{ik}^v)$}}
		-| (9,12) 
		node[pos=.25] {\contour{white}{$hol_k(\gamma^h_{kl})$}}
		node[ pos=.75, rotate=-90] {};
	\node at (10.5,10.5) {$Hol_{k}(\Sigma_k(0))$};
\draw (3,3) -| (9,0)
	node[pos=0.25] {}
	node[pos=0.75, rotate=-90] {\contour{white}{$hol_{l}(\gamma^v_{jl})$}}
		-| (3,3) 
		node[pos=.25] {\contour{white}{$g_{jl}(\gamma^v_{jl}(1))$}}
		node[ pos=.75, rotate=-90] {};	
	\node at (6,1.5) {$Hol_{jl}^{-1}(\gamma^v_{jl})$};
\draw (3,9) -| (9,3)
	node[pos=0.25] {}
	node[pos=0.75, rotate=-90] {\contour{white}{$g_{kl}(x)$}}
		-| (3,9) 
		node[pos=.25] {\contour{white}{$g_{jl}(x)$}}
		node[ pos=.75, rotate=-90] { };	
	\node at (6,6) {$Hol_{ijkl}(x)$};
\draw (3,12) -| (9,9)
	node[pos=0.25] {\contour{white}{$g_{ik}(\gamma^v_{ik}(0))$}}
	node[ pos=0.75,rotate=-90] {\contour{white}{$hol_k(\gamma_{ik}^v)$}}
		-| (3,12) 
		node[pos=.25] {\contour{white}{$g_{ik}(x)$}}
		node[ pos=.75, rotate=-90] {};
	\node at (6,10.5) {$Hol_{ik}^{-1}(\gamma_{ik}^v)$};
\draw (0,3) -| (3,0)
	node[pos=0.25] {}
	node[pos=0.75, rotate=-90] {\contour{white}{$hol_{j}(\gamma_{jl}^v)$}}
		-| (0,3) 
		node[pos=.25] {\contour{white}{$hol_j(\gamma_j^S))$}}
		node[ pos=.75, rotate=-90] {\contour{white}{$hol_{j}(\gamma_j^W)$}};	
	\node at (1.5,1.5) {$Hol_j(\Sigma_j(0))$};
\draw (0,9) -| (3,3)
	node[pos=0.25] {}
	node[pos=0.75, rotate=-90] {\contour{white}{$g_{ij}(x)$}}
		-| (0,9) 
		node[pos=.25] {\contour{white}{$hol_j(\gamma_{ij}^h)$}}
		node[ pos=.75, rotate=-90] {\contour{white}{$g_{ij}(\gamma_{ij}^h(0))$}};	
	\node at (1.5,6) {$Hol_{ij}(\gamma_{ij}^h)$};
\draw (0,12) -| (3,9)
	node[pos=0.25] {\contour{white}{$hol_i(\gamma_{ij}^h)$}}
	node[ pos=0.75,rotate=-90] {\contour{white}{$hol_i(\gamma_{ik}^v)$}}
		-| (0,12) 
		node[pos=.25] {\contour{white}{$hol_i(\gamma_{ij}^h)$}}
		node[ pos=.75, rotate=-90] {\contour{white}{$hol_i(\gamma_i^W)$}};
	\node at (1.5,10.5) {$Hol_i(\Sigma(0))$};
\end{tikzpicture}} = 
\resizebox{6cm}{!}{\begin{tikzpicture}[baseline={([yshift=-.5ex]current bounding box.center)},vertex/.style={anchor=base,
    circle,fill=black!25,minimum size=18pt,inner sep=2pt}]
\draw (9,3) -| (12,0)
	node[pos=0.25] {}
	node[pos=0.75, rotate=-90] {\contour{white}{$hol_{l}(\gamma_l^E)$}}
		-| (9,3) 
		node[pos=.25] {\contour{white}{$hol_l(\gamma_l^S)$}}
		node[ pos=.75, rotate=-90] {};	
	\node at (10.5,1.5) {$Hol_{l}(\Sigma_l(0))$};
\draw (9,9) -| (12,3)
	node[pos=0.25] {}
	node[pos=0.75, rotate=-90] {\contour{white}{$g_{kl}(\gamma^h_{kl}(1))$}}
		-| (9,9) 
		node[pos=.25] {\contour{white}{$hol_l(\gamma^h_{kl})$}}
		node[ pos=.75, rotate=-90] { };	
	\node at (10.5,6) {$Hol_{kl}(\gamma^h_{kl})$};
\draw (9,12) -| (12,9)
	node[pos=0.25] {\contour{white}{$g_{ik}(\gamma_{ik}^v(0))$}}
	node[ pos=0.75,rotate=-90] {\contour{white}{$hol_k(\gamma_{ik}^v)$}}
		-| (9,12) 
		node[pos=.25] {\contour{white}{$hol_k(\gamma^h_{kl})$}}
		node[ pos=.75, rotate=-90] {};
	\node at (10.5,10.5) {$Hol_{k}(\Sigma_k(0))$};
\draw (3,3) -| (9,0)
	node[pos=0.25] {}
	node[pos=0.75, rotate=-90] {\contour{white}{$hol_{l}(\gamma^v_{jl})$}}
		-| (3,3) 
		node[pos=.25] {\contour{white}{$g_{jl}(\gamma^v_{jl}(1))$}}
		node[ pos=.75, rotate=-90] {};	
	\node at (6,1.5) {$Hol_{jl}^{-1}(\gamma^v_{jl})$};
\draw (3,9) -| (9,3)
	node[pos=0.25] {}
	node[pos=0.75, rotate=-90] {\contour{white}{$g_{kl}(x)$}}
		-| (3,9) 
		node[pos=.25] {\contour{white}{$g_{jl}(x)$}}
		node[ pos=.75, rotate=-90] { };	
	\node at (6,6) {$Hol_{ijkl}(y)$};
\draw (3,9) to (4.5, 7.5) to (7.5,7.5) to (9,9);
\draw (3,3) to (4.5, 4.5) to (7.5,4.5) to (9,3);
\draw (4.5, 7.5) to (4.5, 4.5);
\draw (7.5,7.5) to (7.5,4.5);
	\node at (6,8.25) {$Hol_{ik}^{-1}(\rho)$};
	\node[rotate=-90] at (8.25,6) {$Hol_{kl}^{-1}(\rho)$};
	\node at (6,3.75) {$Hol_{jl}(\rho)$};
	\node[rotate=-90]  at (3.76,6) {$Hol_{ij}(\rho)$};
\draw (3,12) -| (9,9)
	node[pos=0.25] {\contour{white}{$g_{ik}(\gamma^v_{ik}(0))$}}
	node[ pos=0.75,rotate=-90] {\contour{white}{$hol_k(\gamma_{ik}^v)$}}
		-| (3,12) 
		node[pos=.25] {\contour{white}{$g_{ik}(x)$}}
		node[ pos=.75, rotate=-90] {};
	\node at (6,10.5) {$Hol_{ik}^{-1}(\gamma_{ik}^v)$};
\draw (0,3) -| (3,0)
	node[pos=0.25] {}
	node[pos=0.75, rotate=-90] {\contour{white}{$hol_{j}(\gamma_{jl}^v)$}}
		-| (0,3) 
		node[pos=.25] {\contour{white}{$hol_j(\gamma_j^S))$}}
		node[ pos=.75, rotate=-90] {\contour{white}{$hol_{j}(\gamma_j^W)$}};	
	\node at (1.5,1.5) {$Hol_j(\Sigma_j(0))$};
\draw (0,9) -| (3,3)
	node[pos=0.25] {}
	node[pos=0.75, rotate=-90] {\contour{white}{$g_{ij}(x)$}}
		-| (0,9) 
		node[pos=.25] {\contour{white}{$hol_j(\gamma_{ij}^h)$}}
		node[ pos=.75, rotate=-90] {\contour{white}{$g_{ij}(\gamma_{ij}^h(0))$}};	
	\node at (1.5,6) {$Hol_{ij}(\gamma_{ij}^h)$};
\draw (0,12) -| (3,9)
	node[pos=0.25] {\contour{white}{$hol_i(\gamma_{ij}^h)$}}
	node[ pos=0.75,rotate=-90] {\contour{white}{$hol_i(\gamma_{ik}^v)$}}
		-| (0,12) 
		node[pos=.25] {\contour{white}{$hol_i(\gamma_{ij}^h)$}}
		node[ pos=.75, rotate=-90] {\contour{white}{$hol_i(\gamma_i^W)$}};
	\node at (1.5,10.5) {$Hol_i(\Sigma_i(0))$};
\end{tikzpicture}}
\end{align}
The above diagram can be ``glued'' together in the following way:
\begin{align}
Hol^{(ijkl)} =& Hol_i(\Sigma_i(0) )\cdot \ldots \cdot \overline{Hol_{ik}^{-1}(\gamma_{ik}^v)} \cdot  \left( \overline{Hol_{ijkl}(x)} \right)\cdot \overline{Hol_{jl}^{-1}(\gamma_{jl}^v} \cdot \ldots \cdot \overline{Hol_l(\Sigma_l(0))} \label{Hol with cube} \\
= &  Hol_i(\Sigma_i(0)) \cdot \ldots \cdot \overline{Hol_{ik}^{-1}(\gamma_{ik}^v)} \nonumber \\
&\cdot  \left( \overline{Hol_{ij}(\rho)} \widehat{Hol_{ik}^{-1}(\rho)} \widehat{Hol_{ijkl}(y)}\overline{Hol_{jl}(\rho)}\widehat{Hol_{kl}^{-1}(\rho)}   \right) \nonumber \\
&\cdot \overline{Hol_{jl}^{-1}(\gamma_{jl}^v}) \cdot \ldots \cdot \overline{Hol_l(\Sigma_l(0))} 
\end{align}
where the notation $\widehat{X}$ is a visual aid for the reader to note that the path-action of $X$, which uses the $r$ direction; i.e. $\widehat{X}$ suggests that the path had to move ``up or down'' in the $r$-direction due to a $Hol(\rho)$ term being placed out of order in the glueing process.   Note that throughout a computation, the $\overline{\bullet}$ and $\hat{\bullet}$ notation is \emph{implicit} and might change when terms are rearranged following the crossed module relations.

The rewrite in \eqref{Hol with cube} demonstrates how the global holonomy is unchanged when the vertex term is replacedwith the remaining $5$ faces of the vertex cube from Proposition \ref{Vertex Cube}.  Now, one can rearrange the terms in \eqref{Hol with cube}, using only the crossed module relation \eqref{eq alpha t h commutes} as follows:
\begin{align}
Hol^{(ijkl)} =  &Hol_i(\Sigma_i(0)) \cdot \overline{Hol_{ij}(\gamma_{ij}^h)}\cdot  \overline{Hol_{ij}(\rho)} \overline{Hol_j(\Sigma_j(0))}\cdot \widehat{Hol_{ik}^{-1}(\gamma_{ik}^v)} \nonumber \\
& \cdot  \cdot \widehat{Hol_{ik}^{-1}(\rho)}  \widehat{Hol_{ijkl}(y)}\cdot \overline{Hol_{jl}(\rho)}\cdot \overline{Hol_{jl}^{-1}(\gamma_{jl}^v) }\\
& \cdot   \cdot \widehat{Hol_k(\Sigma_k(0))} \cdot \widehat{Hol_{kl}^{-1}(\rho)} \cdot \overline{Hol_{kl}(\gamma_{kl}^h)} \cdot \overline{Hol_l(\Sigma_l(0))} \nonumber. \label{eq cube Holijkl}
\end{align}
Here, some faces of the vertex cube were moved to a more convenient\footnote{Convenient in the sense that comparing these non-commutative terms will be easier later on.} location, according to the following order:
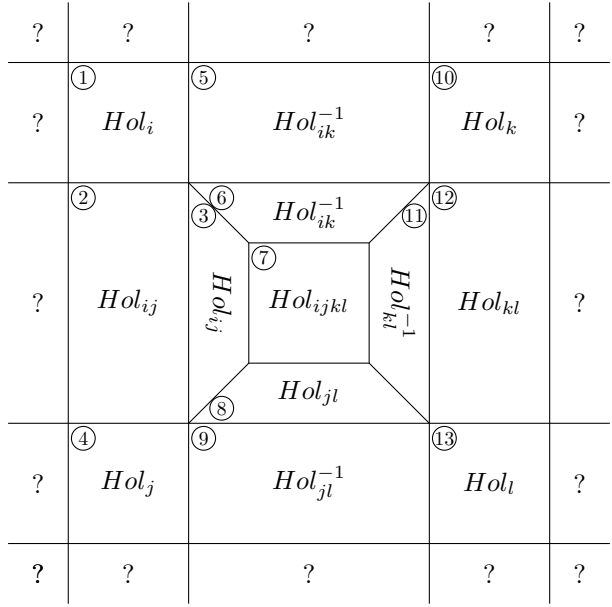
\begin{figure}[H]

\begin{tikzpicture}[scale=0.8]
\foreach \t in {0, 2, 6, 8}
{\draw (\t, -1)--(\t, 9);}
\foreach \s in {0, 2, 6, 8}
{\draw (-1, \s)--(9, \s);}
\draw (2,6)--(3,5)--(5,5)--(6,6);
\draw (2,2)--(3,3)--(5,3)--(6,2);
\draw (3,5)--(3,3);
\draw (5,5)--(5,3);
\node at (1,7) {$Hol_i$};
\node at (1,4) {$Hol_{ij}$};
\node at (1,1) {$Hol_j$};
\node at (4,7) {$Hol_{ik}^{-1}$};
\node at (4,5.5) {$Hol_{ik}^{-1}$};
\node at (4,4) {$Hol_{ijkl}$};
\node at (4,2.5) {$Hol_{jl}$};
\node at (4,1) {$Hol_{jl}^{-1}$};
\node at (7,7) {$Hol_k$};
\node at (7,4) {$Hol_{kl}$};
\node at (7,1) {$Hol_l$};
\node[rotate=-90] at (2.5,4) {$Hol_{ij}$};
\node[rotate=-90] at (5.5,4) {$Hol_{kl}^{-1}$};
 \begin{scope}[every node/.style={scale=.75}]
\node[circle, draw=black, minimum width={width("13")+2pt},inner sep=0pt] at (0.25,7.75) {$1$};
\node[circle, draw=black, minimum width={width("13")+2pt},inner sep=0pt] at (0.25,5.75) {$2$};
\node[circle, draw=black, minimum width={width("13")+2pt},inner sep=0pt] at (0.25,1.75) {$4$};
\node[circle, draw=black, minimum width={width("13")+2pt},inner sep=0pt] at (2.25,7.75) {$5$};
\node[circle, draw=black, minimum width={width("13")+2pt},inner sep=0pt] at (2.25,1.75) {$9$};
\node[circle, draw=black, minimum width={width("13")+2pt},inner sep=0pt] at (6.25,7.75) {$10$};
\node[circle, draw=black, minimum width={width("13")+2pt},inner sep=0pt] at (6.25,5.75) {$12$};
\node[circle, draw=black, minimum width={width("13")+2pt},inner sep=0pt] at (6.25,1.75) {$13$};
\node[circle, draw=black, minimum width={width("13")+2pt},inner sep=0pt] at (2.25,5.45) {$3$};
\node[circle, draw=black, minimum width={width("13")+2pt},inner sep=0pt] at (2.55,5.75) {$6$};
\node[circle, draw=black, minimum width={width("13")+2pt},inner sep=0pt] at (2.55,2.25) {$8$};
\node[circle, draw=black, minimum width={width("13")+2pt},inner sep=0pt] at (5.75,5.45) {$11$};
\node[circle, draw=black, minimum width={width("13")+2pt},inner sep=0pt] at (3.25,4.75) {$7$};
\end{scope};
\foreach \x in {-0.5, 8.5}
{\foreach \y in {-0.5, 1, 4, 7} 
{\node at (\x, \y) {$?$}; \node at(\y, \x) {$?$}; }}
\node at (8.5,8.5) {$?$};
\end{tikzpicture}
\caption{(Vertex Cube) The circled numbers represent the glueing-order of the faces in this arrangement.  The positioning of these numbers points to the ``base point'' of the square from which the target is calculated, as in \eqref{target Holi}.  The question marks represent boundaries or other faces in the larger grid.} \label{Figure VC ordering}
\end{figure}
\noindent Differentiating $Hol^{(\bullet)}$ yields the following lemma:
\begin{lemma}[Vertex Cube Equation, (VCE)]\label{VCE}
For a one-parameter family of squares, $\Sigma_r$, at each vertex $x_{ijkl}$, there is a \emph{vertex cube equation}, denoted by (VCE), where certain terms in the arbitrarily-long product are isolated
\begin{align}
&\ldots \cdot \overline{Hol_{ik}^{-1}} \cdot   \overline{d(Hol_{ijkl})\left( \partiald{r} \right)}  \cdot \overline{Hol_{jl}^{-1}}\cdot \ldots  \\
=&- \ldots \cdot \overline{Hol_{ij}}  \cdot \overline{\alpha_{g_{ij}(x_{ijkl})^{-1}}\left(\restr{a_{ij}}{x_{ijkl}}\left( \partiald{r} \right) \right)} \cdot \overline{Hol_{j}} \ldots  \\
&- \sum\limits_{\bullet_2}\ldots \cdot   \overline{(\alpha_{hol_j \cdot \left(-t(a_{ij})\left( \partiald{r}\right) \right) g_{ij}}\cdot \ldots) (Hol_{\bullet_2})}  \cdot \ldots  \label{term dpath left to Holbullet in VCE} \\
&- \ldots \cdot \overline{Hol_j} \cdot \ldots \cdot   \overline{(\alpha_{hol_j \cdot \left(-t(a_{ij})\left( \partiald{r}\right)\right) g_{ij}}) (Hol_{ik}^{-1})}  \cdot \overline{Hol_{ijkl}} \cdot \ldots \label{term dgij to Holik in VCE} \\
&+  \ldots \cdot \overline{Hol_{ik}^{-1}} \cdot  \overline{\alpha_{g_{ik}^{-1}(x_{ijkl}) }\left(\restr{a_{ik}}{x_{ijkl}}\left( \partiald{r} \right) \right)} \cdot \overline{Hol_{ijkl}} \cdot \ldots   \\
&- \ldots \cdot \overline{Hol_{ik}^{-1}} \cdot \overline{\alpha_{hol_j \cdot \left(A_j\left( \partiald{r} \right) + dg_{ij}\left(\partiald{r} \right)g_{ij}^{-1} \right) \cdot g_{ij}}(Hol_{ijkl})}\cdot \overline{Hol_{jl}^{-1}} \cdot \ldots\\
&-  \ldots \cdot \overline{Hol_{ijkl}} \cdot \overline{ \alpha_{g_{jl}^{-1}(x_{ijkl}) }\left(\restr{a_{jl}}{x_{ijkl}}\left( \partiald{r} \right) \right) } \cdot \overline{Hol_{jl}^{-1}}\cdot \ldots \\
&- \sum\limits_{\bullet_3} \ldots \cdot   \overline{(\alpha_{hol_l \cdot \left(-t(a_{kl})\left( \partiald{r}\right) \right) g_{kl}hol_k\cdot \ldots}) (Hol_{\bullet_3})}  \cdot \ldots \label{term dpath right to Holbullet in VCE}  \\
&-  \ldots \cdot \overline{Hol_{jl}^{-1}}\cdot \overline{(\alpha_{hol_l \cdot \left(-t(a_{kl})\left( \partiald{r}\right) \right) g_{kl}hol_k}) (Hol_{k})}  \cdot \overline{Hol_{kl}} \cdot  \ldots \label{term dpath to Holk in VCE}  \\
&+\ldots \cdot \overline{Hol_k} \cdot \overline{ \alpha_{g_{kl}^{-1}(x_{ijkl}) }\left(\restr{a_{kl}}{x_{ijkl}}\left( \partiald{r} \right) \right)} \cdot \overline{Hol_{kl}} \cdot \ldots  \label{term akl in VCE}.
\end{align}  
where, from Definition \ref{def gerbe}: $$-t(a_{ij}) = A_j + dg_{ij}g_{ij}^{-1} - g_{ij}A_ig_{ij}^{-1}$$ and where $Hol_{\bullet_2}$ is any face in the grid appearing above $Hol_{ik}^{-1}$, and $Hol_{\bullet_3}$ is any face in the grid appearing above $Hol_{k}$.
\end{lemma}
Next, vertical and horizontal edges are considered using the edge cubes from Proposition \ref{Edge Cube}.  The Lemma below yielding the \emph{horizontal} and \emph{vertical edge cube equations} follows from a proof analogous to that of Lemma \ref{VCE} where the following expressions and corresponding diagrams for horizontal edges are used:

For a horizontal edge, $(ij)$, the one parameter family of squares used in this section, $\Sigma(r)$, is restricted to its associated one-parameter family of horizontal edges, $\gamma_{ij}^h(t,r)$.  The edge at height $r$ will be written $\gamma_{ij}^h(-,r)$.  This family of edges also creates a face, $\gamma_{ij}(-, -)$, and two vertical edges, $\gamma_{ij}(0,-)$ and $\gamma_{ij}(1,-)$.  The horizontal portion of Proposition \ref{Edge Cube} yields
\begin{align}\label{h-edge}
Hol =& Hol_{\text{h-edge}}^{(ij)}\\
:=& \ldots \cdot \overline{Hol^{-1}_{pi}}(\gamma_{pi}^v) \cdot \overline{Hol_{pqij}(x_{pqij})} \cdot \overline{Hol_{qj}^{-1}(\gamma_{qj}^v)} \cdot \ldots\\
&\ldots \cdot \overline{Hol_i(\Sigma_i)} \cdot \overline{Hol_{ij}(\gamma_{ij}^h(0,-))} \cdot \widehat{Hol_i(\gamma_{ij}^h(-,-))} \cdot \widehat{Hol_{ij}(\gamma_{ij}^h(-,r))}\\
& \cdot \overline{Hol_{j}^{-1}(\gamma_{ij}^h(-,-))} \cdot \widehat{Hol^{-1}_{ij}(\gamma_{ij}^h(1,-))} \cdot \overline{Hol_j(\Sigma_j)}  \cdot \ldots\\
&\ldots \cdot \overline{Hol_{ik}^{-1}(\gamma_{ik}^v)} \cdot \overline{Hol_{ijkl}(x_{ijkl})} \cdot \overline{Hol_{jl}^{-1}(\gamma_{jl}^v)} \cdot \ldots
\end{align}
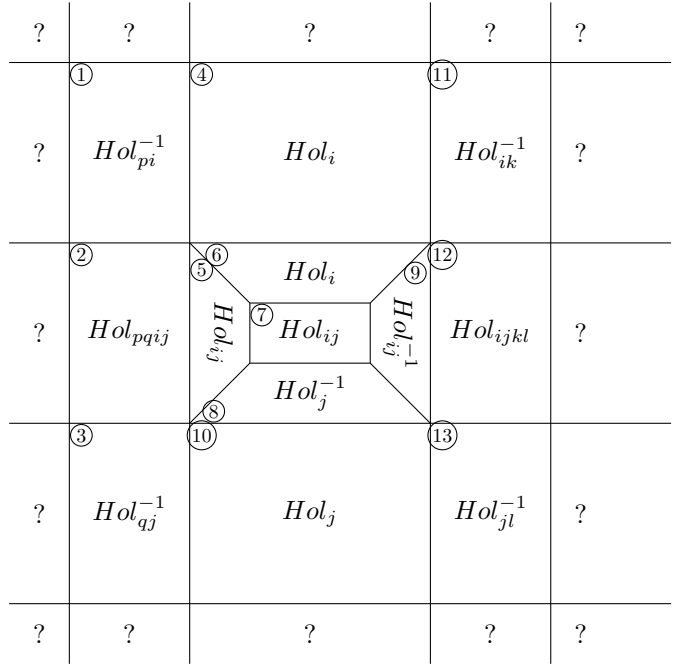
\begin{figure}[H]
\begin{tikzpicture}[scale=0.8]
\foreach \t in {0, 2, 6, 8}
{\draw (\t, -1)--(\t, 10);}
\foreach \s in {0, 3, 6, 9}
{\draw (-1, \s)--(10, \s);}
\draw (2,6)--(3,5)--(5,5)--(6,6);
\draw (2,3)--(3,4)--(5,4)--(6,3);
\draw (3,5)--(3,4);
\draw (5,5)--(5,4);
\node at (1,7.5) {$Hol_{pi}^{-1}$};
\node at (1,4.5) {$Hol_{pqij}$};
\node at (1,1.5) {$Hol_{qj}^{-1}$};
\node at (4,7.5) {$Hol_{i}$};
\node at (4,5.5) {$Hol_{i}$};
\node at (4,4.5) {$Hol_{ij}$};
\node at (4,3.5) {$Hol_{j}^{-1}$};
\node at (4,1.5) {$Hol_{j}$};
\node at (7,7.5) {$Hol_{ik}^{-1}$};
\node at (7,4.5) {$Hol_{ijkl}$};
\node at (7,1.5) {$Hol_{jl}^{-1}$};
\node[rotate=-90] at (2.5,4.5) {$Hol_{ij}$};
\node[rotate=-90] at (5.5,4.5) {$Hol_{ij}^{-1}$};
 \begin{scope}[every node/.style={scale=.75}]
\node[circle, draw=black, inner sep=1pt] at (0.2,8.8) {$1$};
\node[circle, draw=black, inner sep=1pt] at (0.2,5.8) {$2$};
\node[circle, draw=black, inner sep=1pt] at (0.2,2.8) {$3$};
\node[circle, draw=black, inner sep=1pt] at (2.2,8.8) {$4$};
\node[circle, draw=black, inner sep=1pt] at (2.2,2.8) {$10$};
\node[circle, draw=black, inner sep=1pt] at (6.2,8.8) {$11$};
\node[circle, draw=black, inner sep=1pt] at (6.2,5.8) {$12$};
\node[circle, draw=black, inner sep=1pt] at (6.2,2.8) {$13$};
\node[circle, draw=black, inner sep=1pt] at (2.2,5.55) {$5$};
\node[circle, draw=black, inner sep=1pt] at (2.45,5.8) {$6$};
\node[circle, draw=black, inner sep=1pt] at (2.4,3.2) {$8$};
\node[circle, draw=black, inner sep=1pt] at (5.75,5.5) {$9$};
\node[circle, draw=black, inner sep=1pt] at (3.2,4.8) {$7$};
\end{scope};
\foreach \x in {-0.5, 8.5}
{\foreach \y in {-0.5, 1.5, 4.5, 7.5, 9.5} 
{\node at (\x, \y) {$?$}; }}

\foreach \x in {1, 4,7}
{\foreach \y in {-0.5,9.5} 
{\node at (\x, \y) {$?$}; }}
\end{tikzpicture}
\caption{(horizontal edge cube)} \label{Figure hEC ordering}
\end{figure}

Similarly, there is a one parameter family of vertical edges, $\gamma_{ij}^v(-,-)$, and thus the following expressions and diagram for vertical edges: 
\begin{align}\label{hol v edge}
Hol = Hol_{\text{v-edge}}^{(ij)}:= &\ldots \cdot \overline{Hol_{pi}(\gamma_{pi}^h}) \cdot \overline{Hol^{-1}_i(\gamma_{ij}^v(-,0))} \cdot \overline{Hol_i(\Sigma_i)} \cdot \overline{Hol_{ik}(\gamma_{ik}^h}) \cdot \ldots\\
&\ldots \cdot \widehat{Hol_{piqj}}(x_{pqij}) \cdot \widehat{Hol_{ij}^{-1}}(\gamma_{ij}^v(1,-)) \cdot \widehat{Hol_{ij}^{-1}}(\gamma_{ij}^v(-,r)) \cdot \overline{Hol_{ij}} (\gamma_{ij}^h(0,-))\cdot \overline{Hol_{ikjl}(x_{ikjl})} \cdot \ldots\\
&\ldots \cdot \widehat{Hol_{qj}(\gamma_{qj}^h)} \cdot \widehat{Hol_{j}(\Sigma_j)} \cdot \widehat{Hol_{j}(\gamma_{ij}^v(-,r)} \cdot \overline{Hol_{jl}(\gamma_{jl}^h)} \cdot \ldots
\end{align}
\begin{figure}[H]\label{v edge}
\begin{tikzpicture}[scale=0.8]
\foreach \s in {0, 2, 6, 8}
{\draw (-1, \s)--(10, \s);}
\foreach \t in {0, 3, 6, 9}
{\draw (\t, -1)--(\t, 9);}
\draw (3,6)--(4,5)--(5,5)--(6,6);
\draw (3,2)--(4,3)--(5,3)--(6,2);
\draw (4,5)--(4,3);
\draw (5,5)--(5,3);
\node at (1.5,7) {$Hol_{pi}$};
\node at (1.5,4) {$Hol_{i}$};
\node at (1.5,1) {$Hol_{ik}$};
\node at (4.5,7) {$Hol_{pqij}$};
\node at (4.5,5.5) {$Hol_{ij}^{-1}$};
\node[rotate=-90] at (4.5,4) {$Hol_{ij}^{-1}$};
\node at (4.5,2.5) {$Hol_{ij}$};
\node at (4.5,1) {$Hol_{ikjl}$};
\node at (7.5,7) {$Hol_{qj}$};
\node at (7.5,4) {$Hol_{j}$};
\node at (7.5,1) {$Hol_{jl}$};
\node[rotate=-90] at (3.5,4) {$Hol_{i}^{-1}$};
\node[rotate=-90] at (5.5,4) {$Hol_{j}$};
 \begin{scope}[every node/.style={scale=.75}]
\node[circle, draw=black, inner sep=1pt] at (0.2,7.8) {$1$};
\node[circle, draw=black, inner sep=1pt] at (0.2,5.8) {$3$};
\node[circle, draw=black, inner sep=1pt] at (0.2,1.8) {$4$};
\node[circle, draw=black, inner sep=1pt] at (3.2,7.8) {$5$};
\node[circle, draw=black, inner sep=1pt] at (3.2,1.8) {$9$};
\node[circle, draw=black, inner sep=1pt] at (6.2,7.8) {$10$};
\node[circle, draw=black, inner sep=1pt] at (6.2,5.8) {$11$};
\node[circle, draw=black, inner sep=1pt] at (6.2,1.8) {$13$};
\node[circle, draw=black, inner sep=1pt] at (3.2,5.55) {$2$};
\node[circle, draw=black, inner sep=1pt] at (3.45,5.8) {$6$};
\node[circle, draw=black, inner sep=1pt] at (3.4,2.2) {$8$};
\node[circle, draw=black, inner sep=1pt] at (5.75,5.5) {$12$};
\node[circle, draw=black, inner sep=1pt] at (4.2,4.8) {$7$};
\end{scope};
\foreach \x in {-0.5, 9.5}
{\foreach \y in {-0.5, 1, 4.5, 7, 8.5} 
{\node at (\x, \y) {$?$}; }}
\foreach \x in {1.5, 4.5, 7.5}
{\foreach \y in {-0.5, 8.5} 
{\node at (\x, \y) {$?$}; }}
\end{tikzpicture}
\caption{(vertical edge cube)} \label{Figure vEC ordering}
\end{figure}
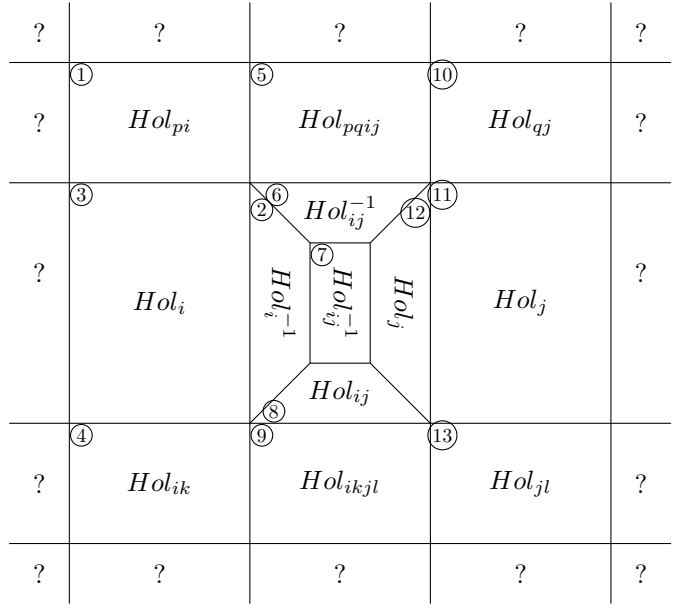
Differentiating these equations gives the following lemma:
\begin{lemma}[Edge Cube Equations, (ECEh) and (ECEv)]\label{ECE}
For a one-parameter family of squares, $\Sigma_r$, the following \emph{horizontal edge cube equation} at $\gamma_{ij}^h$, denoted by (ECEh), holds where certain terms in the arbitrarily-long product are isolated
\begin{align}
&\ldots \cdot \overline{Hol_{i}} \cdot   \overline{d(Hol_{ij})\left( \partiald{r} \right)}  \cdot \overline{Hol_{j}}\cdot \ldots \nonumber \\
=&- \ldots \cdot \overline{Hol_{i}}  \cdot \overline{ \alpha_{g_{ij}(\gamma_{ij}^h(0))^{-1}}\left(\restr{a_{ij}}{\gamma_{ij}^h(0)} \left(\partiald{r} \right) \right) } \cdot \overline{Hol_{ij}} \cdot \ldots \label{term west aij in ECEh} \\
&- \ldots \overline{Hol_{i}}  \cdot \overline{\int_{\gamma_{ij}^h}\alpha_*(B_i)\left(\partiald{r}\right)  } \cdot \overline{Hol_{ij}}  \cdot \ldots  \\
&- \ldots  \cdot \overline{Hol_{i}}  \cdot \overline{(\alpha_{\left(\restr{A_j}{\gamma_{ij}^h(0)} \left(\partiald{r} \right) + \restr{dg_{ij}}{\gamma_{ij}^h(0)}\left(\partiald{r} \right) g_{ij}^{-1} \right)g_{ij}})(Hol_{ij})} \cdot \overline{Hol_{j}}  \cdot  \ldots \nonumber  \\
&+ \ldots \cdot \overline{Hol_{ij}}  \cdot \overline{\int_{\gamma_{ij}^h}\alpha_*(B_j)\left(\partiald{r}\right)  } \cdot \overline{Hol_{j}}  \cdot  \ldots  \nonumber \\
&+ \ldots \cdot \overline{Hol_{ij}}  \cdot \overline{\alpha_{g_{ij}(\gamma_{ij}^h(1))^{-1}}\left(\restr{a_{ij}}{\gamma_{ij}^h(1)} \left(\partiald{r} \right) \right)} \cdot \overline{Hol_{j}}  \cdot  \ldots . \nonumber
\end{align}
Similarly, the following \emph{vertical edge cube equation} at $\gamma_{ij}^v$, denoted by (ECEv), holds:
\begin{align}
&\ldots \cdot  \overline{Hol_{piqj}} \cdot   \overline{d(Hol_{ij}^{-1}) \left( \partiald{r} \right)}  \cdot \overline{Hol_{ikjl}}\cdot \ldots \\
=&+ \ldots \cdot \overline{Hol_{pi}}   \cdot \overline{\int_{\gamma_{ij}^v}\alpha_*(B_i)\left(\partiald{r}\right)  }  \cdot \overline{Hol_{i}} \cdot \ldots \label{term Bi in ECEv} \\
&-\sum\limits_{\bullet_2} \ldots \cdot \overline{(\alpha_{\left(\restr{\widetilde{dhol_i}}{\gamma_{ij}^v}\left( \partiald{r} \right) \right)hol_i\cdot g_{pi}\cdot \ldots})(Hol_{\bullet_2})} \cdot  \ldots \label{term dpath to Holbullet2 in ECEv} \\
&- \ldots \cdot \overline{Hol_{ik}}   \cdot \ldots \cdot \overline{(\alpha_{\left(\restr{\widetilde{dhol_i}}{\gamma_{ij}^v}\left( \partiald{r} \right)  \right)hol_i\cdot g_{pi}})(Hol_{piqj})} \cdot \overline{Hol_{ij}^{-1}} \cdot \ldots \label{term dholi to Holpiqj in ECEv} \\
&+ \ldots \cdot \overline{Hol_{piqj}} \cdot  \overline{ \alpha_{g_{ij}(\gamma_{ij}^v(0))^{-1}}\left(\restr{a_{ij}}{\gamma_{ij}^v(0)} \left(\partiald{r} \right) \right) }   \cdot \overline{Hol_{ij}^{-1}}  \cdot  \ldots  \\
&- \ldots \cdot  \overline{Hol_{piqj}}\cdot    \overline{(\alpha_{\left(\restr{A_i}{\gamma_{ij}^v(1)} \left(\partiald{r} \right) + \restr{dhol_i}{\gamma_{ij}^v}\left(\partiald{r} \right) hol_i^{-1}  \right)})(Hol_{ij}^{-1})} \cdot \overline{Hol_{ikjl}} \cdot \ldots  \\
&- \ldots \cdot \overline{Hol^{-1}_{ij}} \cdot \overline{ \alpha_{g_{ij}(\gamma_{ij}^v(1))^{-1}}\left(\restr{a_{ij}}{\gamma_{ij}^v(1)} \left(\partiald{r} \right) \right) }   \cdot \overline{Hol_{ikjl}}  \cdot  \ldots  \\
&-\sum\limits_{\bullet_3} \ldots \cdot   \overline{(\alpha_{g_{jl} \left(\restr{\widetilde{dhol_j}}{\gamma_{ij}^v}\left( \partiald{r} \right)  \right)hol_j\cdot g_{qj}\cdot \ldots})(Hol_{\bullet_3})}  \cdot  \ldots \label{term dpath to Holbullet3 in ECEv} \\
&- \ldots \cdot \overline{Hol_{ikjl}}   \cdot \ldots \cdot   \overline{(\alpha_{g_{jl} \left(\restr{\widetilde{dhol_j}}{\gamma_{ij}^v}\left( \partiald{r} \right)  \right)hol_j\cdot g_{qj}})(Hol_{qj})}  \cdot \overline{Hol_{j}} \cdot  \ldots \label{term dhol to Holqj in ECEv} \\
&- \ldots \cdot \overline{Hol_{qj}}  \cdot \overline{(\alpha_{g_{jl} \left(\restr{\widetilde{dhol_j}}{\gamma_{ij}^v}\left( \partiald{r} \right) \right)hol_j})(Hol_{j})} \cdot \overline{Hol_{jl}}\cdot \ldots \label{term dhol to Holj in ECEv}  \\
&- \ldots \cdot \overline{Hol_{j}} \cdot  \overline{\int_{\gamma_{ij}^v}\alpha_*(B_j)\left(\partiald{r}\right)}    \cdot \overline{Hol_{jl}} \cdot \ldots.
\end{align}
where the abbreviation $\tilde{d hol}$ is given by  $$\restr{\widetilde{dhol}}{\gamma}\left( \partiald{r} \right):= \restr{A}{\gamma(1)} \left(\partiald{r} \right) + \restr{dhol}{\gamma}\left(\partiald{r} \right) hol^{-1}  - hol \restr{A}{\gamma(0)}\left(\partiald{r} \right) hol^{-1}$$
and where $Hol_{\bullet_2}$ is any square in the grid appearing above $Hol_{piqj}$, and $Hol_{\bullet_3}$ is any square in the grid appearing above $Hol_{qj}$,
\end{lemma}

\subsection{Proof of Theorem \ref{dhol nice}}\label{dHol nice proof}
In this section, it is shown that $d(Hol)$ does not have differential information at the interior edges or vertices. At the same time, the boundary terms are gathered into a convenient expression.  It is important to note that $d(Hol)$, applied to $\restr{\partiald{r}}{r=0}$, can be written, using the Liebniz rule as in \eqref{Liebniz of Hol}.  This expression has the following types of terms
\begin{enumerate}
\item \emph{3-curvature terms}: namely the terms in $\overline{d(Hol_i)}$ from Lemma \ref{local lemma} where one $B_i$ is replaced with an $H_i$
\item \emph{Boundary terms}: these are terms which occur on the boundary of $\Sigma$.
\item \emph{Edge and Vertex terms}: \label{edge and vertex terms}
\begin{enumerate}
\item The remaining 2 types of terms coming from $\overline{d(Hol_i)}$: four \emph{edge terms}, written as one integral around the four sides in Lemma \ref{local lemma}, and a \emph{corner term} where $A_i$ is applied to the upper left corner of the square.
\item $\overline{d(Hol_{ij})}$ and $\overline{d(Hol_{ijkl})}$: These are the de Rham differentials applied to any (vertical or horizontal) edges or vertices, see Lemmas \ref{VCE} and \ref{ECE}.
\item \emph{Path terms}: These are the terms which appear from differentiating the $\alpha_{\text{path}}$ of any $\overline{Hol_{\bullet}}= \alpha_{\text{path}}(Hol_{\bullet})$ term in $d(Hol)$ (see Lemma \ref{lemma d alpha}).
\end{enumerate}
\end{enumerate}
\subsubsection{Summary of Cancelation}
For the moment, all of the terms which end up canceling are collected and compared.  Begin by using Lemmas \ref{VCE} and \ref{ECE} to replace any $\restr{\partiald{r}}{r=0} Hol_{ij}$ or $\restr{\partiald{r}}{r=0} Hol_{ijkl}$ with the corresponding expression and then showing that the $\widehat{Hol_{ij}}$ (\emph{edge}) terms cancel with each other.  In other words, for each term $\overline{Hol_{\bullet}}$ or $\widehat{Hol_{\bullet}}$, applying $\partiald{r}$ to such terms and using Lemma \ref{lemma d alpha}, yields terms of the form $$\overline{\restr{\partiald{r}}{r=0} Hol_{\bullet}}.$$  The following table gives a summary of all instances where two such terms will appear in $d(Hol)$ with opposite sign.  Note that for the reader's convenience  the \emph{Local Lemma} \ref{local lemma} is labelled as ``$(LL)$'', the \emph{Vertex Cube Equation} from Lemma \ref{VCE} as ``$(VCE)$'', and the \emph{vertical} and \emph{horizontal edge cube equations} from Lemma \ref{ECE} as ``$(ECEv)$'' and ``$(ECEh)$'', respectively.  
\begin{table}[H]\label{Interior side terms table}
\begin{tabular}{l | l | l}
 \hline
     Label \quad 
    & Term  
    & Found in\\
    \hline
 \multicolumn{3}{c}{Interior \emph{side} terms} \\
 \hline
    (A1) 
    & $ \ldots \overline{Hol_p} \cdot \overline{Hol_{pi} }\cdot \overline{\int_{\gamma_{pi}^h} \alpha_* B_i} \cdot \overline{Hol_i} \ldots $ \dotfill
    & (LL) (ECEh)\\
    (A2) 
    & $ \ldots \overline{Hol_i} \cdot \overline{\int_{\gamma_{ij}^h} \alpha_* B_i}\cdot \overline{Hol_{ij}} \cdot \overline{Hol_j} \ldots $ \dotfill
    & (LL) (ECEh)\\
    (A3) 
    & $\ldots \overline{Hol_i} \cdot \overline{\int_{\gamma^v_{iq} }\alpha_* B_i}\cdot \ldots \cdot \overline{Hol^{-1}_{iq}} \cdot \ldots \cdot \overline{Hol_q} $ \dotfill
    & (LL) (ECEv)\\  
    (A4) 
    & $\ldots \overline{Hol_k} \cdot \ldots \cdot \overline{Hol^{-1}_{ki}}\cdot \ldots \cdot \overline{\int_{\gamma^v_{ki}} \alpha_* B_i} \cdot \overline{Hol_i} \ldots $ \dotfill
    & (LL) (ECEv)\\   
    (B1) 
    & $ \ldots \overline{Hol_p} \cdot \ldots \cdot \overline{\restr{a_{pi}}{\gamma_{pi}^v(0)}} \cdot  \overline{Hol_{pi}}  \cdot \ldots \cdot \overline{Hol_i} \ldots $ \dotfill
    & (VCE) (ECEv)\\ 
    (B2)  
    & $ \ldots \overline{Hol_p} \cdot \ldots \cdot  \overline{Hol_{pi}} \cdot \overline{\restr{a_{pi}}{\gamma_{pi}^v(1)}}  \cdot \ldots \cdot \overline{Hol_i} \ldots $ \dotfill
    & (VCE) (ECEv)\\  
    (B3)  
    & $ \ldots \overline{Hol_i} \cdot \overline{\restr{a_{ij}}{\gamma_{ij}^h(0)}} \cdot \overline{Hol_{ij}}   \cdot \overline{Hol_j} \ldots $ \dotfill
    & (VCE) (ECEh)\\  
    (B4)  
    & $ \ldots \overline{Hol_i}  \cdot \overline{Hol_{ij}} \cdot \overline{\restr{a_{ij}}{\gamma_{ij}^h(1)}}  \cdot \overline{Hol_j} \ldots $ \dotfill
    & (VCE) (ECEh)\\      
 \hline
\end{tabular}    
\end{table}
The other terms which appear from applying $\restr{\partiald{r}}{r=0}$ to each term $\overline{Hol_{\bullet}}$ or $\widehat{Hol_{\bullet}}$ are the \emph{path terms} from Lemma \ref{lemma d alpha}, 
$$\alpha_{\left( \restr{\partiald{r}}{r=0}\dots  \right) } (\overline{Hol_{\bullet}}),$$
and are labelled ``path-$d(Hol)$''.

The two tables below summarize all of the instances where $\emph{path terms}$ show up an even number of times, with opposite signs.  In the first table, the focus is on \emph{path terms} resulting from $d(Hol)$; i.e. they show up one time in differentiating the path approaching each term in the expression \eqref{eq global Hol} for $Hol^{\mathcal{N}}$.  The second table deals with \emph{path terms} which cancel amongst the other relations.
\begin{table}[H]\label{table Interior path terms}
\begin{tabular}{l | l | l}
 \hline
     Label \quad 
    & Term  
    & Found in\\
    \hline
 \multicolumn{3}{c}{Interior \emph{path} terms} \\
 \hline
    (C1) 
    & $\ldots \cdot \alpha_{(\ldots \cdot d hol_i \cdot \ldots )}(Hol_z) \cdot \ldots $ \dotfill
    & (path-$d(Hol)$) (ECEv)\\
    (C2) 
    & $\ldots \cdot \alpha_{(\ldots \cdot d hol_i \cdot \ldots )}(Hol_{yz}) \cdot \ldots $ \dotfill
    & (path-$d(Hol)$) (ECEv)\\
    (C3) 
    & $\ldots \cdot \alpha_{(\ldots \cdot d hol_i \cdot \ldots )}(Hol^{-1}_{wx}) \cdot \ldots $ \dotfill
    & (path-$d(Hol)$) (ECEv)\\ 
    (C4) 
    & $\ldots \cdot \alpha_{(\ldots \cdot d hol_i \cdot \ldots )}(Hol_{wxyz}) \cdot \ldots $ \dotfill
    & (path-$d(Hol)$) (ECEv)\\   
    (D1) 
    & $\ldots \cdot \alpha_{(\ldots \cdot d g_{ij} \cdot \ldots )}(Hol_z) \cdot \ldots $ \dotfill
    & (path-$d(Hol)$) (VCE)\\
    (D2) 
    & $\ldots \cdot \alpha_{(\ldots \cdot d g_{ij} \cdot \ldots )}(Hol_{yz}) \cdot \ldots $ \dotfill
    & (path-$d(Hol)$) (VCE)\\
    (D3) 
    & $\ldots \cdot \alpha_{(\ldots \cdot d g_{ij} \cdot \ldots )}(Hol^{-1}_{yx}) \cdot \ldots $ \dotfill
    & (path-$d(Hol)$) (VCE)\\ 
    (D4) 
    & $\ldots \cdot \alpha_{(\ldots \cdot d g_{ij}  \cdot \ldots )}(Hol_{wxyz}) \cdot \ldots $ \dotfill
    & (path-$d(Hol)$) (VCE)\\        
 \hline
\end{tabular}    
\end{table}
Note that in (C3) one can have $wx = ix$ and in (D4) one can have $ij = wx$.  In (D2) one can have $ij = yz$ but then the cancellation is due to (path-$d(Hol)$) and (ECEh).   
\begin{table}[H]\label{table Interior Corner terms}

\begin{tabular}{l | l | l}
 \hline
     Label \quad 
    & Term  
    & Found in\\
    \hline
 \multicolumn{3}{c}{Interior \emph{corner} terms} \\
 \hline
    (E1) 
    & $\ldots \cdot \alpha_{(\ldots \cdot hol_i(\gamma^v_{qi}) \cdot \restr{A_i}{\gamma^v_{qi}(0)} \cdot \ldots )}(Hol_z) \cdot \ldots $ \dotfill
    & (ECEv) (VCE) \\
    (E2) 
    & $\ldots \cdot \alpha_{(\ldots \cdot hol_i(\gamma^v_{qi}) \cdot \restr{A_i}{\gamma^v_{qi}(0)} \cdot \ldots )}(Hol_{yz}) \cdot \ldots $ \dotfill
    & (ECEv) (VCE) \\   
    (E3) 
    & $\ldots \cdot \alpha_{(\ldots \cdot  hol_i(\gamma^v_{iq}) \restr{A_i} {\gamma^v_{iq}(0)} \cdot \ldots )}(Hol_{wx}^{-1}) \cdot \ldots $ \dotfill
    & (ECEv) (VCE)\\ 
    (E4) 
    & $\ldots \cdot \alpha_{(\ldots \cdot  hol_i(\gamma^v_{iq}) \restr{A_i} {\gamma^v_{iq}(0)}  \cdot \ldots )}(Hol_{wxyz}) \cdot \ldots $ \dotfill
    & (ECEv) (VCE)\\     
    (F1) 
    & $\ldots \cdot \alpha_{(\ldots \cdot \restr{A_i} {\gamma^v_{qi}(1)} hol_i(\gamma^v_{qi}) \cdot \ldots )}(Hol_z) \cdot \ldots $ \dotfill
    & (ECEv) (VCE)\\ 
    (F2) 
    & $\ldots \cdot \alpha_{(\ldots \restr{A_i} {\gamma^v_{qi}(1)}  hol_i(\gamma^v_{qi}) \cdot \ldots )}(Hol_{yz}) \cdot \ldots $ \dotfill
    & (ECEv) (VCE)\\ 
    (F3) 
    & $\ldots \cdot \alpha_{(\ldots \cdot \restr{A_i} {\gamma^v_{iq}(1)}  hol_i (\gamma^v_{iq})\cdot \ldots )}(Hol^{-1}_{wx}) \cdot \ldots $ \dotfill
    & (ECEv) (VCE)\\ 
    (F4) 
    & $\ldots \cdot \alpha_{(\ldots \cdot \restr{A_i} {\gamma^v_{iq}(1)}  hol_i(\gamma^v_{iq}) \cdot \ldots )}(Hol_{wxyz}) \cdot \ldots $ \dotfill
    & (ECEv) (VCE)\\ 
    (G1) 
    & $\ldots \cdot \alpha_{(\ldots \cdot  hol_i(\gamma^v_{qi}) \restr{A_i} {\gamma^v_{qi}(0)} )}(Hol_{i}) \cdot \ldots $ \dotfill
    & (ECEv) (LL)\\          
    (G2) 
    & $\ldots \cdot \alpha_{(\ldots \cdot   \restr{A_i} {\gamma^v_{pi}(1)} \cdot g_{pi})}(Hol_{pi}) \cdot \ldots $ \dotfill
    & (ECEv) (ECEh)\\  
 \hline
\end{tabular}    
\end{table}
Note that in (E3) one can have $wx=iq$, in (E4) one can have $wxyz=wiyq$, and in (F1) one can have $z=i$.

The only terms that are left, for the interior of $\Sigma$, after all of this cancelation are the $3$-curvature terms, $H_i$ integrated over the square; i.e. the term $\int H$ in Theorem \ref{dhol nice}.

\subsubsection{Getting the Boundary Right}
Modulo the boundary, it has thus far been shown that $d(Hol) \equiv Hol \cdot \int_{Sq} H$.  In order to get the boundary terms to work out properly, it remains to check that all terms that accumulate at the boundary of $\Sigma$, due to not being able to cancel with a missing adjacent square, either cancel or are of the form $B_i$ or $a_{ij}$ as described in the statement of Theorem \ref{dhol nice}.  In particular, consider again the general grid from Definition \ref{def of Hol}.  Note that when $d(Hol_{ij})$ is replaced with its corresponding cube equation at the northern boundary of the square, $\Sigma$, a $Hol_{ij}$ is differentiated as it collapses to the northern boundary, placing an $a_{ij}$ at that spot.  In such cases, it will be useful to write any term 
\begin{equation}\label{a term northern boundary}
Hol_a \cdot \ldots \cdot \overline{Hol_m} \cdot \overline{a_{ab}\left(\partiald{r}\right)} \cdot  \overline{Hol_{ab}^{-1}} \cdot \ldots \cdot \overline{Hol_p}
\end{equation} as $\overline{a_{ab}} \cdot Hol$, simply by changing the path-action, $\overline{\bullet}$, for the $a_{ij}$ term.  The easier terms to deal with will be on the Northern and Eastern boundaries.  By using the algebra of the crossed module, the equation \eqref{a term northern boundary} can be rewritten as desired.  However, explaining this algebra is a lot easier by recalling $hh' = \alpha_{t(h)}(h') h$ and observing the equality as a picture:
\begin{equation}
 \resizebox{6cm}{!}{\begin{tikzpicture}[scale=1.2, baseline={([yshift=-.5ex]current bounding box.center)},vertex/.style={anchor=base,
     circle,fill=black!25,minimum size=18pt,inner sep=2pt}]   
\def\x{3}
 \path[gray]  (0,\x)       edge (6,\x);
  \path[gray]  (6,\x)       edge (6,0);
   \path[gray]  (6,0)       edge (0,0);
    \path[gray]  (0,0)       edge (0,\x);
  \path[gray]  (3,\x)       edge (3,0);
\node[gray]  at (1.5, \x/2) {$Hol_a \cdot \ldots \cdot \overline{Hol_m}$};
\node[gray]  at (4.5, \x/2) {$\overline{Hol_{ab}^{-1}} \cdot \ldots \cdot \overline{Hol_p}$};
\node[above right,]  at (3, \x) {$a_{ab}$};
\node[draw=gray,above, circle]  at (1.5, \x/2 ) {\color{gray}{$1$}};
\node[draw=gray, above, circle]  at (4.5, \x/2) {\color{gray}{$3$}};
\node[draw=black, above left,yshift=5pt,  circle]  at (3, \x) {$2$};
\path[thick] (0,\x) edge (0,0);
\path[thick] (0,0) edge (3,0);
\path[thick] (3,0) edge (3,\x);
\end{tikzpicture}}
=  \resizebox{6cm}{!}{\begin{tikzpicture}[scale=1.2, baseline={([yshift=-.5ex]current bounding box.center)},vertex/.style={anchor=base,
     circle,fill=black!25,minimum size=18pt,inner sep=2pt}]   
     \def\x{3}
 \path[gray]  (0,\x)       edge (6,\x);
  \path[gray]  (6,\x)       edge (6,0);
   \path[gray]  (6,0)       edge (0,0);
    \path[gray]  (0,0)       edge (0,\x);
  \path[gray]  (3,\x)       edge (3,0);
\node[gray]  at (1.5, \x/2) {$Hol_a \cdot \ldots \cdot \overline{Hol_m}$};
\node[gray]  at (4.5, \x/2) {$\overline{Hol_{ab}^{-1}} \cdot \ldots \cdot \overline{Hol_p}$};
\node[above right,]  at (3, \x) {$a_{ab}$};
\node[draw=gray,above, circle]  at (1.5, \x/2) {\color{gray}{$2$}};
\node[draw=gray, above, circle]  at (4.5, \x/2) {\color{gray}{$3$}};
\node[draw=black, above left, yshift=5pt, circle]  at (3, \x) {$1$};
\path[thick] (0,\x) edge (3,\x);
\end{tikzpicture}}
\end{equation}   
In a similar fashion, based on the ordering of the Edge Cube equations (Lemmas \ref{VCE} and \ref{ECE}) and the Local Lemma \ref{local lemma}, all of the $B_i$ and $a_{ij}$ terms appearing along the Northern and Eastern boundary can be factored outside of $Hol$.

For the Western and Southern boundaries, there is one extra tool needed.  Considering again an $a_{ij}$ term, let us consider the term
\begin{equation}\label{a term western boundary}
Hol_a \cdot \overline{Hol_{ae}} \cdot \overline{Hol_e}\cdot  \overline{Hol_{ei}} \cdot \overline{a_{ei}\left(\partiald{r}\right)}\cdot \overline{Hol_i}  \cdot \ldots \cdot \overline{Hol_p}
\end{equation} 
which will be useful to rewrite as $Hol \cdot \overline{a_{ei}}$.  The tool here is to realize there are left-over terms on the boundary which assemble precisely to ``$-[a_{ei}, - ]$''.  To see this, one last type of term coming from $d(Hol)$ is finally used, which has not been previously incorporated: the derivative of the path-action terms along the Western and Southern boundaries coming from each $\overline{\bullet}$ in the expression for $Hol$.  Note that these terms did not appear for the Northern and Eastern boundaries since the convention is to use the path approaching a term going along the Western boundary, then along the Southern boundary, and then up towards the desired location through the interior.  In other words, \eqref{eq global Hol} can be written as
\begin{align}
&Hol_a \cdot \overline{Hol_{ae}} \cdot \overline{Hol_e} \cdot  \overline{Hol_{ei}} \cdot  \cdot{Hol_i} \cdot \ldots \cdot \overline{Hol_p}\\
= &Hol_a \cdot \ldots \overline{Hol_{ei} } \cdot \alpha_{hol_a^{-1} \cdot hol_{ae}^{-1} \cdot hol_e^{-1} \cdot g_{ei}^{-1} }(Hol_i) \cdot \ldots \cdot \overline{Hol_p}
\intertext{which will momentarily be written as}
= &Hol_1 \cdot \overline{Hol_{ei} } \cdot \alpha_{hol_a^{-1} \cdot hol_{ae}^{-1} \cdot hol_e^{-1} \cdot g_{ei}^{-1} }(Hol_i) \cdot  \overline{Hol_2}.
\end{align} 
Using all of the various terms occurring at this Southern $ei$ boundary-corner, they can be combined in a useful way, where the reference to which set of relations the term comes from is listed in place of an equation label:
\begin{align}
&-Hol_1 \cdot \overline{Hol_{ei}} \cdot \alpha_{(\ldots \cdot g_{ei}^{-1}(dg_{ei}g_{ei}^{-1}))}(Hol_i \cdot \overline{Hol_2}) \tag{$d(Hol)$-path}\\
&-Hol_1 \cdot \alpha_{(\ldots g_{ei}^{-1}(A_i + dg_{ei}g_{ei}^{-1}))}(\alpha_{g_{ei}}(Hol_{ei})) \cdot \overline{Hol_i} \cdot \overline{Hol_2} \tag{ECEh} \\
&-Hol_1 \overline{Hol_{ei}} \cdot \alpha_{(\ldots \cdot g_{ei}^{-1}\cdot A_i)}(Hol_i) \cdot \overline{Hol_2} \tag{LL}\\
&-Hol_1 \cdot \alpha_{g_{ei}^{-1}}(a_{ei}) \cdot \overline{Hol_{ei}} \cdot \overline{Hol_i} \cdot \overline{Hol_2} \tag{ECEh}\\
&+Hol_1 \cdot \alpha_{(\ldots \cdot A_e)}(Hol_{ei} \cdot \overline{Hol_i} \cdot \overline{Hol_2}) \tag{$d(Hol)$-path}\\
&-Hol_1 \cdot \overline{Hol_{ei}} \cdot \overline{Hol_i} \cdot \alpha_{(\ldots \cdot g_{ei}^{-1}A_i hol_i)}(Hol_2) \tag{$d(Hol)$-path}\\
=&-Hol_1 \cdot \alpha_{g_{ei}^{-1}}(a_{ei}) \cdot \overline{Hol_{ei}} \cdot \overline{Hol_i} \cdot \overline{Hol_2}\nonumber \\
&+ Hol_1 \cdot \alpha_{(\ldots \cdot(A_e - dg_{ei}g_{ei}^{-1} - g_{ei}^{-1}A_ig_{ei}))}(\alpha_{g_{ei}}(Hol_{ei}) \cdot \overline{Hol_j} \cdot \overline{Hol_2}) \nonumber \\
=&-Hol_1 \cdot \alpha_{g_{ei}^{-1}}(a_{ei}) \cdot \overline{Hol_{ei}} \cdot \overline{Hol_i} \cdot \overline{Hol_2}\nonumber  \\
&+ Hol_1 \cdot \alpha_{(\ldots \cdot g_{ei}^{-1}(t(a_{ei})))}(\alpha_{g_{ei}}(Hol_{ei}) \cdot \overline{Hol_j} \cdot \overline{Hol_2})\nonumber  \\
=&-Hol_1 \cdot \alpha_{g_{ei}^{-1}}(a_{ei}) \cdot \overline{Hol_{ei}} \cdot \overline{Hol_i} \cdot \overline{Hol_2}\nonumber  \\
&+ Hol_1 \cdot \alpha_{(\ldots \cdot g_{ei}^{-1})}\left[ a_{ei}, \alpha_{g_{ei}}(Hol_{ei}) \cdot \overline{Hol_j} \cdot \overline{Hol_2}\right]\nonumber \\
= &Hol \cdot \overline{a_{ei}}\nonumber 
\end{align}
This calculation above demonstrates the appearance of the term, 
$$Hol \cdot \left(\sum\limits_{\partial Sq} a \right),$$
from equation \eqref{dhol nice}.

To demonstrate the appearance of the term, 
$$Hol \cdot \left(\int_{\partial Sq} B   \right),$$
from equation \eqref{dhol nice}, a similar technique is applied to the $B_i$ integrated along the Western and Southern boundaries using the vanishing fake curvature condition $t(B_i) = R_i$, which is now shown below.

In a similar manner to the above, first write \eqref{eq global Hol} as
\begin{equation*} Hol_1 \cdot \overline{Hol_n} \cdot \alpha_{(\ldots \cdot hol_n^{-1})}(\overline{Hol_2}). \end{equation*}
After differentiating, \emph{side} terms along the path, $\gamma_n^{S}$, are obtained which can be rewritten in the desirable fashion:
\begin{align}
&-Hol_1 \cdot \overline{Hol_n} \cdot \overline{\int_{\gamma_n^{S}} \alpha_*(B_n)} \cdot \overline{Hol_2}\\
&+Hol_1 \cdot \overline{Hol_n} \cdot \alpha_{\ldots \cdot \int_{\gamma_n^{S}} \alpha_*(R_n)} ( \overline{Hol_2}) \nonumber \\
=&-Hol_1 \cdot \overline{Hol_n} \cdot \overline{\int_{\gamma_n^{S}} \alpha_*(B_n)} \cdot \overline{Hol_2} \nonumber \\
&+Hol_1 \cdot \overline{Hol_n} \cdot \alpha_{\ldots \cdot t\left(\int_{\gamma_n^{S}} \alpha_*(t(B_n))\right)} ( \overline{Hol_2}) \nonumber \\
=&-Hol_1 \cdot \overline{Hol_n} \cdot \overline{\int_{\gamma_n^{S}} \alpha_*(B_n)} \cdot \overline{Hol_2} \nonumber \\
&+Hol_1 \cdot \overline{Hol_n} \cdot \left( \overline{\left[\int_{\gamma_n^{S}} \alpha_*(B_n) ,Hol_2 \right]}\right)\nonumber \\
= & Hol \cdot \overline{\int_{\gamma_n^{S}} \alpha_*(B_n)} \nonumber
\end{align}

\section{Special Cases}
The results of this paper are now briefly stated as three special cases: 
\begin{enumerate}
\item The surface holonomy of spheres, $M^{S^2} \to H$.
\item The surface holonomy, $M^{Sq} \to H$, which uses a crossed module $(H \xrightarrow{t} G, \alpha)$, whose $\alpha$-action is given by inner-automorphisms.
\item  The surface holonomy for abelian gerbes.
\end{enumerate}

\subsection{Surface Holonomy of Spheres}\label{sec spheres}
First some useful propositions for $Hol^{\mathcal{N}}$ for the case when $\Sigma \in M^{S^2}$ are recorded below.  By way of applying Proposition \ref{prop t of Hol and composition} to the case where the eastern and western boundaries are identified with each other, the northern edge is collapsed to a point, and the southern edge is collapsed to a point, we obtain the following two corollaries just as one can find in \cite{MP2} and \cite{SWIII}:
\begin{cor}\label{hol on sphere is in center}
The 2-holonomy of a sphere takes its values in the center of the Lie group, $H$.
\end{cor}

\begin{cor}\label{Hol transforms}
The transformation of 2-holonomy between open sets in $M^{S^2}$ is given by
$$Hol^{\mathcal{N}_{I'}} (\Sigma) = \alpha_{g_{ii'}(0,0)}(Hol^{\mathcal{N}_I})$$
for $\Sigma \in \mathcal{N}_{I} \cap \mathcal{N}_{I'} \subset M^{S^2}$.
\end{cor}
Similarly, in the case of $M^{S^2}$ Theorem \ref{dhol nice} can be simplified to an original corollary as follows.
\begin{cor}\label{dhol nice for spheres}
The total derivative of 2-holonomy, $d(Hol)$, can be written:
\begin{equation}  d(Hol)=  -(\alpha_{Hol})_*(A_{i_{(1,1)}}) +Hol\cdot \int_{S^2} H \end{equation}
where $Hol=Hol^\mathcal {N}:\mathcal{N}\to H\subset Mat$ is defined on the open set $\mathcal N\subset M^{S^2}$.
\end{cor}
Note that $\int_{S^2} H \in \Omega^1(\mathcal{N} , \mathfrak{h})$ is a one-form on the open subset $\mathcal{N} \subset M^{S^2}$; where $\mathfrak{h}$ is the Lie algebra of $H$, as usual.  For two open subsets $\mathcal{N}_I, \mathcal{N}_{I'} \subset M^{S^2}$, the transformation of $\int_{S^2} H$ in $\mathcal{N}_I \cap \mathcal{N}_{I'} \subset M^{S^2}$ can be written as follows.
\begin{cor}
The integral of the 3-curvature over a sphere transforms in the following way:
$$\int_{S^2} H^{\mathcal{N}_I}  = \alpha_{g_{ii'}(0,0)}(\int_{S^2} H^{\mathcal{N}_{I'}})$$
where by $H^{\mathcal{N}_{I'}}$ we mean to use the local 3-curvature as defined by the local data on $\mathcal{N}_{I'}$.
\end{cor}
\subsubsection{$\alpha : G \to Inn(H)$}\label{section alpha inner}
In this section the following special case is assumed, which will allow for considerable simplifications.
\begin{setting}\label{set alpha inner}
Suppose that $\alpha$ factors through the inner automorphisms of $H$: 
\begin{equation}
\xymatrix{
G \ar[dr]_{\alpha} \ar[r]^{\alpha} &Aut(H)\\
&Inn(H) \ar@{^{(}->}[u]}
\end{equation}
\end{setting}
\begin{remark}\label{simple setting}
Note that if $Y \in Z(\mathfrak{h})$ then in Setting \ref{set alpha inner}, $\alpha_g (Y) = Y$ for any $g \in G$.  Similarly, if $h \in Z(H)$, then $\alpha_X (h) = 0$ for any $X \in \mathfrak{g}$, and for $Y \in Z(\mathfrak{h})$, $X \in \mathfrak{g}$ then $\alpha_{X}(Y) = 0$.  
\end{remark}
Recall from Proposition \ref{hol on sphere is in center} that the $2$-holonomy, a function on $\mathcal{N} \subset M^{S^2}$, has trivial target yielding $Hol^{\mathcal{N}} \in Z(H)$, where $Z(H)$ is the center of $H$.  It was shown in Proposition \ref{Hol transforms} that $Hol$ transforms between open subsets of $M^{S^2}$ via $\alpha_{g_{ij}}$ but since the domain is in $M^{S^2}$ and the action factors through an inner automorphism, it follows that $Hol^{\mathcal{N}}$ functions agree on overlaps (also proven in \cite{MP2}):
\begin{cor}
The function $Hol: M^{S^2} \to H$ given by $\restr{Hol}{\mathcal{N}}:= Hol^{\mathcal{N}}$ is well-defined (globally). 
\end{cor}

\begin{cor}\label{Cor H forms glue}
In the case of Setting \ref{set alpha inner}, the local $3$-curvature forms, $H_i$, glue together to a global 3-form, $H \in \Omega^3(M, \mathfrak{h})$.
\end{cor}

As a consequence of the above fact we obtain the final original result of this paper,
\begin{cor}\label{prop dHol global sphere}
The total derivative of 2-holonomy, $d(Hol)= Hol \cdot \int_{S^2} H \in \Omega^1(M^{S^2}, Mat)$, is globally defined.\end{cor}

\subsection{Abelian Gerbes and Surface Holonomy}\label{sec abelian}
For this section, assume that squares, $\Sigma: [0,1]^2 \to M$, are closed surfaces in $M$.  Recall that in an abelian gerbe, the structure crossed module $(H \to 1, t)$ is used, where $H$ is an abelian group, and the trivial group, $1$, takes the place of the group, $G$.  Note that in this case the $\alpha$-action is trivial: $\alpha(h) = h$, and so in the notation of this paper, it follows that $\overline{Hol_{\bullet}}= Hol_{\bullet}$.  Such a gerbe, $\mathcal{G}$, has $g_{ij} =1$, $A_i = 0$.  Note, then that we can recover the analogous results found in \cite{TWZ}:
\begin{cor}
The surface holonomy function is well defined on $M$. In particular we have $Hol^{\mathcal{N}_{I'}} = Hol^{\mathcal{N}_{I}}$.
\end{cor}

\begin{cor}
For the function $Hol=: M \to H\subset Mat$ defined on any open set $\mathcal N\subset M^{Sq}$, we can compute its total derivative $d(Hol)\in \Omega^1(M^{Sq},Mat)$ as follows:
\begin{equation}
d(Hol)= Hol\cdot \int_{Sq} H
\end{equation}
where $H$ is the global $3$-form of Corollary \ref{Cor H forms glue}.
\end{cor}


\end{document}